\documentclass[10pt]{amsart}
\usepackage{latexsym,amsfonts,amsmath,amssymb}
\newtheorem{theorem}{Theorem}
\newtheorem{corollary}[theorem]{Corollary}
\newtheorem{sublemma}{Lemma}[theorem]
\newtheorem{lemma}[theorem]{Lemma}
\newtheorem{question}[theorem]{Question}
\newtheorem{observation}[theorem]{Observation}
\newtheorem{claim}[theorem]{Claim}
\newtheorem{subclaim}{Claim}[sublemma]
\newtheorem{conjecture}[theorem]{Conjecture}
\newtheorem{fact}[theorem]{Fact}
\newtheorem{definition}[theorem]{Definition}
\newtheorem{remark}[theorem]{Remark}
\newtheorem{example}[theorem]{Example}
\newtheorem{exercise}{Exercise}[section]
\def\Theorem #1.#2 #3\par{\def\claimname{Theorem}\setbox1=\hbox{#1}\ifdim\wd1=0pt
   \begin{theorem}{\rm #2} #3\end{theorem}\else
   \newtheorem{#1}[theorem]{#1}\begin{#1}\label{#1}{\rm #2} #3\end{#1}\fi}
\def\Corollary #1.#2 #3\par{\def\claimname{Corollary}\setbox1=\hbox{#1}\ifdim\wd1=0pt
   \begin{corollary}{\rm #2} #3\end{corollary}\else
   \newtheorem{#1}[theorem]{#1}\begin{#1}\label{#1}{\rm #2} #3\end{#1}\fi}
\def\Lemma #1.#2 #3\par{\def\claimname{Lemma}\setbox1=\hbox{#1}\ifdim\wd1=0pt
   \begin{lemma}{\rm #2} #3\end{lemma}\else
   \newtheorem{#1}[theorem]{#1}\begin{#1}\label{#1}{\rm #2} #3\end{#1}\fi}
\def\SubLemma #1.#2 #3\par{\def\claimname{Lemma}\setbox1=\hbox{#1}\ifdim\wd1=0pt
   \begin{sublemma}{\rm #2} #3\end{sublemma}\else
   \newtheorem{#1}{#1}[theorem]\begin{#1}\label{#1}{\rm #2} #3\end{#1}\fi}
\def\Question #1.#2 #3\par{\def\claimname{Question}\setbox1=\hbox{#1}\ifdim\wd1=0pt
   \begin{question}{\rm #2} #3\end{question}\else
   \newtheorem{#1}[theorem]{#1}\begin{#1}\label{#1}{\rm #2} #3\end{#1}\fi}
\def\Observation #1.#2 #3\par{\def\claimname{Observation}\setbox1=\hbox{#1}\ifdim\wd1=0pt
   \begin{observation}{\rm #2} #3\end{observation}\else
   \newtheorem{#1}[theorem]{#1}\begin{#1}\label{#1}{\rm #2} #3\end{#1}\fi}
\def\Claim #1.#2 #3\par{\def\claimname{Claim}\setbox1=\hbox{#1}\ifdim\wd1=0pt
   \begin{claim}{\rm #2} #3\end{claim}\else
   \newtheorem{#1}[theorem]{#1}\begin{#1}\label{#1}{\rm #2} #3\end{#1}\fi}
\def\SubClaim #1.#2 #3\par{\def\claimname{Claim}\setbox1=\hbox{#1}\ifdim\wd1=0pt
   \begin{subclaim}{\rm #2} #3\end{subclaim}\else
   \newtheorem{#1}{#1}[sublemma]\begin{#1}\label{#1}{\rm #2} #3\end{#1}\fi}
\def\Conjecture #1.#2 #3\par{\def\claimname{Conjecture}\setbox1=\hbox{#1}\ifdim\wd1=0pt
   \begin{conjecture}{\rm #2} #3\end{conjecture}\else
   \newtheorem{#1}[theorem]{#1}\begin{#1}\label{#1}{\rm #2} #3\end{#1}\fi}
\def\Fact #1.#2 #3\par{\def\claimname{Fact}\setbox1=\hbox{#1}\ifdim\wd1=0pt
   \begin{fact}{\rm #2} #3\end{fact}\else
   \newtheorem{#1}[theorem]{#1}\begin{#1}\label{#1}{\rm #2} #3\end{#1}\fi}
\def\Definition #1.#2 #3\par{\def\claimname{Definition}\setbox1=\hbox{#1}\ifdim\wd1=0pt
   \begin{definition}{\rm #2} {\rm #3}\end{definition}\else
   \newtheorem{#1}[theorem]{#1}\begin{#1}\label{#1}{\rm #2} {\rm #3}\end{#1}\fi}
\def\Remark #1.#2 #3\par{\def\claimname{Remark}\setbox1=\hbox{#1}\ifdim\wd1=0pt
   \begin{remark}{\rm #2} {\rm #3}\end{remark}\else
   \newtheorem{#1}[theorem]{#1}\begin{#1}\label{#1}{\rm #2} {\rm #3}\end{#1}\fi}
\def\Example #1.#2 #3\par{\def\claimname{Example}\setbox1=\hbox{#1}\ifdim\wd1=0pt
   \begin{example}{\rm #2} #3\end{example}\else
   \newtheorem{#1}[theorem]{#1}\begin{#1}\label{#1}{\rm #2} #3\end{#1}\fi}
\def\Exercise #1.#2 #3\par{\def\claimname{Exercise}\setbox1=\hbox{#1}\ifdim\wd1=0pt
   {\footnotesize\begin{exercise}{\rm #2} {\rm #3}\end{exercise}}\else
   \newtheorem{#1}[section]{#1}{\footnotesize\begin{#1}\label{#1}{\rm #2} {\rm #3}\end{#1}}\fi}
\def\QuietTheorem #1.#2 #3\par{\setbox1=\hbox{#1}\ifdim\wd1=0pt\proclaim{Theorem {\rm #2}}{#3}\else\proclaim{#1 {\rm #2}}{#3}\fi}
\newcommand{\proclaim}[2]{\smallskip\noindent{\bf #1} {\sl#2}\par\smallskip}
\def\Proclaim #1.#2 #3\par{\proclaim{#1 {\rm #2}}{#3}}
\newif\ifinproof\inprooffalse
\def\Proof#1. {\setbox1=\hbox{#1}\ifdim\wd1=0pt\begin{proof}\else\begin{proof}[#1]\fi\inprooftrue} 
\newcommand{\QED}{\end{proof}}

\def\BF#1.{{\bf #1.}}
\newcommand{\cal}{\mathcal}
\def\Abstract #1\par{\begin{abstract}#1\end{abstract}}
\def\Title #1\par{\title{#1}\maketitle}
\def\Author #1\par{\author{#1}}
\def\Acknowledgement#1\par{\thanks{#1}}
\def\Chapter #1\par{\chapter{#1}}
\def\Section #1\par{\section{#1}}
\def\QuietSection #1\par{\section*{#1}}
\def\SubSection #1\par{\subsection{#1}}
\def\SubSubSection #1\par{\subsubsection{#1}}
\def\MidTitle #1\par{\bigskip\goodbreak\centerline{\small\bf #1}\bigskip\noindent}
\def\Margin #1\par{\marginpar{\tiny #1}}

\newcommand{\citindex}{\index[citdep]}
\newcommand{\Ref}[2][Theorem]{\edef\entry{\ref{#2}..@#1 \ref{#2}!\ifinproof\arabic{chapter}.\arabic{section}.\arabic{theorem}..@\claimname\ \arabic{chapter}.\arabic{section}.\arabic{theorem}\else \arabic{chapter}.\arabic{section}..@Section \arabic{chapter}.\arabic{section}\fi|dotfill\ p.}\expandafter\citindex\expandafter{\entry}\ref{#2}}
\def\bottomnote #1\par{{\renewcommand{\thefootnote}{}\footnotetext{#1}}}
\newcommand{\F}{{\mathbb F}}

\newcommand{\N}{{\mathbb N}}

\newcommand{\Q}{{\mathbb Q}}
\newcommand{\Z}{{\mathbb Z}}
\def\R{{\mathbb R}}

\newcommand{\id}{\mathop{\hbox{\small id}}}

\newfont{\msam}{msam10 at 12pt}
\newcommand{\from}{\mathbin{\vbox{\baselineskip=2pt\lineskiplimit=0pt
                         \hbox{.}\hbox{.}\hbox{.}}}}
\newcommand{\of}{\subseteq}

\newcommand{\set}[1]{\{\,{#1}\,\}}
\newcommand{\singleton}[1]{\left\{{#1}\right\}}

\newcommand{\elesub}{\prec}

\newcommand{\inverse}{{-1}}
\newcommand{\jump}{{\!\triangledown}}

\newcommand{\dom}{\mathop{\rm dom}}

\newcommand{\ran}{\mathop{\rm ran}}

\newcommand{\Add}{\mathop{\rm Add}}

\newcommand{\Th}{\mathop{\rm Th}}

\newcommand{\image}{\mathbin{\hbox{\tt\char'42}}}

\newcommand{\restrict}{\upharpoonright}
\newcommand{\satisfies}{\models}

\newcommand{\proves}{\vdash}

\newcommand{\converges}{\downarrow}

\newcommand{\union}{\cup}

\newcommand{\Union}{\bigcup}

\newcommand{\intersect}{\cap}

\newcommand{\trianglelt}{\lhd}

\newcommand{\smalllt}{\mathrel{\mathchoice{\raise2pt\hbox{$\scriptstyle<$}}{\raise1pt\hbox{$\scriptstyle<$}}{\raise0pt\hbox{$\scriptscriptstyle<$}}{\scriptscriptstyle<}}}
\newcommand{\smallleq}{\mathrel{\mathchoice{\raise2pt\hbox{$\scriptstyle\leq$}}{\raise1pt\hbox{$\scriptstyle\leq$}}{\raise1pt\hbox{$\scriptscriptstyle\leq$}}{\scriptscriptstyle\leq}}}

\newcommand{\boolval}[1]{\mathopen{\lbrack\!\lbrack}\,#1\,\mathclose{\rbrack\!\rbrack}}
\def\[#1]{\boolval{#1}}
\newcommand{\gcode}[1]{\ulcorner\!#1\!\urcorner}
\newcommand{\UnderTilde}[1]{{\setbox1=\hbox{$#1$}\baselineskip=0pt\vtop{\hbox{$#1$}\hbox to\wd1{\hfil$\sim$\hfil}}}{}}
\newcommand{\Undertilde}[1]{{\setbox1=\hbox{$#1$}\baselineskip=0pt\vtop{\hbox{$#1$}\hbox to\wd1{\hfil$\scriptstyle\sim$\hfil}}}{}}
\newcommand{\undertilde}[1]{{\setbox1=\hbox{$#1$}\baselineskip=0pt\vtop{\hbox{$#1$}\hbox to\wd1{\hfil$\scriptscriptstyle\sim$\hfil}}}{}}
\newcommand{\UnderdTilde}[1]{{\setbox1=\hbox{$#1$}\baselineskip=0pt\vtop{\hbox{$#1$}\hbox to\wd1{\hfil$\approx$\hfil}}}{}}
\newcommand{\Underdtilde}[1]{{\setbox1=\hbox{$#1$}\baselineskip=0pt\vtop{\hbox{$#1$}\hbox to\wd1{\hfil\scriptsize$\approx$\hfil}}}{}}

\newcommand{\st}{\mid}
\renewcommand{\th}{{\hbox{\scriptsize th}}}

\newcommand{\minus}{\setminus}
\newcommand{\iso}{\cong}
\def\<#1>{\langle#1\rangle}

\newcommand{\TC}{\mathop{\hbox{\sc tc}}}

\newcommand{\WO}{\mathop{\hbox{\sc wo}}}
\newcommand{\HC}{\mathop{\hbox{\sc hc}}}
\newcommand{\LC}{\mathop{\hbox{\sc lc}}}
\newcommand{\ZFC}{\hbox{\sc zfc}}

\newcommand{\CH}{\hbox{\sc ch}}

\newcommand{\GCH}{\hbox{\sc gch}}

\newcommand{\PA}{\hbox{\sc pa}}
\newcommand{\TA}{\hbox{\sc ta}}

\newcommand{\cell}[1]{\boxit{\hbox to 17pt{\strut\hfil$#1$\hfil}}}
\newcommand{\head}[2]{\lower2pt\vbox{\hbox{\strut\footnotesize\it\hskip3pt#2}\boxit{\cell#1}}}
\newcommand{\boxit}[1]{\setbox4=\hbox{\kern2pt#1\kern2pt}\hbox{\vrule\vbox{\hrule\kern2pt\box4\kern2pt\hrule}\vrule}}
\newcommand{\Col}[3]{\hbox{\vbox{\baselineskip=0pt\parskip=0pt\cell#1\cell#2\cell#3}}}
\newcommand{\tapenames}{\raise 5pt\vbox to .7in{\hbox to .8in{\it\hfill input: \strut}\vfill\hbox to
.8in{\it\hfill scratch: \strut}\vfill\hbox to .8in{\it\hfill output: \strut}}}
\newcommand{\Head}[4]{\lower2pt\vbox{\hbox to25pt{\strut\footnotesize\it\hfill#4\hfill}\boxit{\Col#1#2#3}}}
\newcommand{\Dots}{\raise 5pt\vbox to .7in{\hbox{\ $\cdots$\strut}\vfill\hbox{\ $\cdots$\strut}\vfill\hbox{\
$\cdots$\strut}}}
\renewcommand{\dots}{\raise5pt\hbox{\ $\cdots$}}
\newcommand{\factordiagramup}[6]{$$\begin{array}{ccc}
#1&\raise3pt\vbox{\hbox to60pt{\hfill$\scriptstyle
#2$\hfill}\vskip-6pt\hbox{$\vector(4,0){60}$}}&#3\\ \vbox
to30pt{}&\raise22pt\vtop{\hbox{$\vector(4,-3){60}$}\vskip-22pt\hbox
to60pt{\hfill$\scriptstyle #4\qquad$\hfill}}
     &\ \ \lower22pt\hbox{$\vector(0,3){45}$}\ {\scriptstyle #5}\\
\vbox to15pt{}&&#6\\
\end{array}$$}
\newcommand{\factordiagram}[6]{$$\begin{array}{ccc}
#1&&\\ \ \ \raise22pt\hbox{$\vector(0,-3){45}$}\ {\scriptstyle #2}
&\raise22pt\hbox{$\vector(2,-1){90}$}\raise5pt\llap{$\scriptstyle#3$\qquad\quad}&\vbox
to25pt{}\\ #4&\raise3pt\vbox{\hbox to90pt{\hfill$\scriptstyle
#5$\hfill}\vskip-6pt\hbox{$\vector(4,0){90}$}}&#6\\
\end{array}$$}
\newcommand{\df}{\it} 
\newcommand{\RL}{\mathcal{R}}
\newcommand{\Rw}{\mathcal{R}_w}
\newcommand{\Rc}{\mathcal{R}_f}
\newcommand{\Rf}{\mathcal{R}_f}
\newcommand{\Ra}{\mathcal{R}_a}
\newcommand{\Rev}{\mathcal{R}_e}
\newtheorem{proposition}[theorem]{Proposition}

\begin{document}
\author{Joel David Hamkins}
\address{J. D. Hamkins, The College of Staten Island of The City  University of New York, Mathematics, 2800 Victory Boulevard, Staten  Island, NY 10314
\& The Graduate Center of The City University of New York, Ph.D.~Program in Mathematics, 365 Fifth Avenue, New York, NY 10016, USA}
\email{jhamkins@gc.cuny.edu, http://jdh.hamkins.org}

\author{Russell Miller}
\address{R. Miller, Queens College of The City University of New York,  Mathematics, 65-30 Kissena Boulevard, Flushing, New York 11367
\& The Graduate Center of The City University of New York, Ph.D.~Program in Computer Science, 365 Fifth Avenue, New York, NY 10016, USA}
\email{rmiller@forbin.qc.edu}

\author{Daniel Seabold}
\address{D. Seabold, Department of Mathematics, Hofstra University,  Hempstead, NY 11549-1030, USA} \email{matdes@hofstra.edu}

\author{Steve Warner}
\address{S. Warner, Department of Mathematics, Hofstra University,  Hempstead, NY 11549-1030, USA} \email{matsjw@hofstra.edu}

\bottomnote MSC: 03D60; 03D45; 03C57; 03E15. Keywords: infinite time  Turing machines, computable model theory. The research of the first two
authors has been supported in part by grants from the Research Foundation of CUNY,  and the first author is additionally thankful to the Institute
for Logic, Language and Computation and the NWO (Bezoekersbeurs B62-612) for  supporting his summer 2005 stay at Universiteit van Amsterdam.

\Abstract We introduce infinite time computable model theory, the computable model theory arising with infinite time Turing machines,  which provide
infinitary notions of computability for structures built on the reals  $\R$. Much of the finite time theory generalizes to the infinite time context,
but several fundamental questions, including the infinite time  computable analogue of the Completeness Theorem, turn out to be  independent of \ZFC.

\Title Infinite time computable model theory

\Section Introduction

\markboth{J. D. HAMKINS, R. MILLER, D. SEABOLD AND S. WARNER}{INFINITE  TIME COMPUTABLE MODEL THEORY} Computable model theory is model theory with a
view to the  computability of the structures and theories that arise (for a standard reference, see \cite{HandbookOfRecMath1998}). Infinite time
computable model theory, which we introduce here, carries out this program with the infinitary notions of computability provided by infinite time
Turing machines. The motivation for a broader context is that, while finite time computable model theory is necessarily limited to countable models
and theories, the infinitary context naturally allows for uncountable models and theories, while retaining the computational nature of the
undertaking. Many constructions generalize from finite time computable model theory, with  structures built on $\N$, to the infinitary theory, with
structures built on $\R$. In this article, we introduce the basic theory and consider the infinitary  analogues of the completeness theorem, the
L\"owenheim-Skolem Theorem, Myhill's theorem and others. It turns out that, when stated in their fully general infinitary forms, several of these
fundamental questions are independent of \ZFC. The analysis makes use of techniques both from computability theory and  set theory. This article
follows up \cite{Hamkins2005:InfinitaryComputabilityWithITTM}.

\SubSection Infinite time Turing machines

The definitive introduction to infinite time Turing machines appears in  \cite{HamkinsLewis2000:InfiniteTimeTM}, but let us quickly describe how they
work. The $$\tapenames\Head100{start}\Col100\Col100\Col000\Col100\Col000\Dots$$ hardware of an infinite time Turing machine is  identical to a
classical (three tape) Turing machine, with a head  reading and writing $0$s and $1$s on the one-way infinite tapes, following the instructions  of a
finite program with finitely many states. Computation begins with  the input on the {\it input} tape and the head on the left-most cell in the {\it
start} state. Successor steps of computation are determined by the  program in exactly the classical manner. At any limit ordinal stage, as a matter
of definition, the machine resets the head to the left-most cell, assumes the {\it limit} state and updates the tape so that every cell exhibits the
$\limsup$ of the previous  values displayed in that cell. This is equivalent to using the limit value, if the value displayed by the cell has
stabilized, and otherwise $1$. Computation ceases only when the {\it halt} state is explicitly obtained, and in this case the output is whatever is
written on the output tape. (If the head falls off the tape, no output is given.) If $p$ is a program, it computes a function $\varphi_p$, defined by
$\varphi_p(x)=y$ if and only if on  input $x$ the computation determined by $p$ leads to output $y$. The natural context here for input and output is
the Cantor space ${}^\omega 2$ of  all infinite binary sequences, which we will denote by $\R$ and refer to as the set of reals. A (partial) function
$f\from\R\to\R$ is infinite time {\df computable}  if it is $\varphi_p$ for some program $p$. Binary and $n$-ary functions can be equivalently
modelled either by adding additional input tapes, or by  viewing a single real as the interleaving of the digits of $n$ many  reals. A set $A\of\R$
is infinite time {\df decidable} if its characteristic  function is infinite time computable. The set $A$ is infinite time {\df semi-decidable} if
the function $1\restrict A$ with domain $A$ and constant value $1$ is computable. In this article, we will freely use the terms {\df computable} and
{\df decidable} to mean infinite time computable and infinite time decidable, though we will sometimes specify ``infinite time'' for  clarity. When
referring to the classical notions of computability, we  will always say ``finite time computable'' and ``finite time decidable.'' We regard  the
natural numbers $\N$ as coded in $\R$ by identifying $n$ with the binary sequence consisting of $n$ ones followed by zeros. A real is {\df writable}
if it is $\varphi_p(0)$ for some program $p$. A real is {\df accidentally} writable if it appears on one of the tapes during any computation
$\varphi_p(0)$. A real is {\df eventually} writable if it appears on the output tape of a (not necessarily halting) computation $\varphi_p(0)$, and
from some point on in that computation, it is never changed.  An ordinal $\alpha$ is {\df clockable} if there is a computation $\varphi_p(0)$ moving
to the {\it  halt} state exactly on the $\alpha^\th$ computational step.

The growing body of literature on infinite time Turing machines includes \cite{HamkinsLewis2000:InfiniteTimeTM}, \cite{Welch2000:LengthsOfITTM},
\cite{Welch2000:Eventually}, \cite{HamkinsSeabold2001:OneTape}, \cite{Loewe2001:RevisionSequences}, \cite{HamkinsLewis2002:PostProblem},
\cite{Hamkins2002:Turing}, \cite{HamkinsWelch2003:Pf=NPf}, \cite{DeolalikarHamkinsSchindler2005:NPcoNP},
\cite{Hamkins2005:InfinitaryComputabilityWithITTM}, \cite{Welch2005:ActionOfOneTapeMachines}, \cite{Koepke2005:TuringComputationsOnOrdinals}.

\SubSection Basic definitions

The main idea will be that a computable model is one whose underlying  set is decidable and whose functions and relations are uniformly  computable.
In order to make this precise, let us first be more specific about our syntax  and how it is represented. A language consists of a collection of
function, relation and constant symbols, with each function and relation symbol  assigned a finite arity. In addition, every language has the logical
connective symbols $\wedge$, $\vee$, $\neg$, $\to$, $\leftrightarrow$, parentheses, the equality symbol $=$, quantifiers $\forall$, $\exists$,
variable symbols $v_0$, $v_1$, and so on. In finite time computable model theory, in  order to bring these syntactic objects into the realm of
computability,  one views each symbol in the (countable) language as being represented by a  particular natural number, its G\"odel code, so that the
various  syntactic objects---such as terms, formulas and sentences---are simply finite  sequences of these codes, which can in turn be coded with a
single  natural number.

Infinite time computable model theory, however, offers the possibility  of {\it uncountable} computable models. And because we will want to  consider
the elementary or atomic diagrams of such models, the possibility of  uncountable languages is unavoidable. Clearly, we cannot expect to code  such
languages using G\"odel codes only in $\N$. Therefore, we work in a more general  context, where the symbols of a language are represented with
G\"odel  codes in $\R$, rather than $\N$. This conforms with the philosophy of infinite  time computability, where the fundamental inputs and outputs
of  computations are real numbers. A {\df computable presentation} of a language $\mathcal  L$ is the assignment of a G\"odel code $\gcode{s}$ to
every function, relation and constant symbol $s$ in the language, in such a way that the set of  such codes for symbols in $\mathcal L$ is decidable,
and there are  computable functions telling us, given any $\gcode{s}$, what kind of symbol $s$ is  and, when it is a function or relation symbol,
what arity it has. We  assume that the basic logical symbols (logical connectives, $=$ symbol, parentheses, variable symbols, quantifiers) have
simple G\"odel codes in $\N$.

Given the G\"odel codes of the underlying symbols, one develops the  G\"odel coding of all the usual syntactic notions. For example, a term  $\tau$
is a particular kind of finite sequence of function, constant and variable  symbols, and we may assign the G\"odel code $\gcode{\tau}$ via the  usual
manner of coding finite sequences of reals with reals. Similarly, any formula  $\varphi$ in the language is a finite sequence of symbols from the
language, and we can assign it a natural G\"odel code. We assume that the G\"odel  coding of the language is undertaken in such a way that we can
unambiguously determine whether a given G\"odel code is the code of a formula or an individual symbol, and what kind; that from the G\"odel code of a
formula or term we can compute the G\"odel codes of the subformulas and  subterms; and that the G\"odel codes are uniquely readable. For any
computable presentation $\mathcal L$, it follows that all the elementary syntactic  notions are computable from the G\"odel codes, such as finding
the  inductive construction history of a formula or term or determining whether a given occurrence of a variable in a formula is free or not.

\Definition. In the infinite time context, a {\df computable model} is  a structure ${\mathcal A}=\<A,f^{\mathcal A},R^{\mathcal A},c^{\mathcal
A}>_{f,R,c\in\mathcal L}$ in a language $\cal L$, with a fixed  computable presentation of $\cal L$, such that the underlying set  $A\of\R$ of the
model is decidable and the functions, relations and constants of $\mathcal A$  are uniformly computable from their input and the G\"odel codes of
their symbols. A structure has a {\df computable presentation} if it is  isomorphic to a computable model.

A simple recursive argument shows that the value of any term $\tau(\vec  a)$ is uniformly computable from its G\"odel code $\gcode{\tau}$ and  the
input $\vec a$. It follows that one can compute the truth in $\cal A$ of any given atomic  formula. Specifically, the atomic diagram of $\cal A$ is
the set $\Delta_0({\mathcal A})=\set{\varphi[\vec a]\st\varphi\hbox{ atomic},\vec a\in A^{<\omega}, {\mathcal  A}\satisfies\varphi[\vec a]}$, and if
$\mathcal A$ is a computable  model, then we can decide, on input $\gcode{\varphi}$ and $\vec a$, whether  $\varphi[\vec a]\in\Delta_0({\mathcal
A})$. More generally, we define:

\Definition. A model $\mathcal A$ is {\df (infinite time) decidable} if  the full elementary diagram of the structure $\Delta({\mathcal
A})=\set{\varphi[\vec a]\st {\mathcal A}\satisfies\varphi[\vec a]}$ is infinite time decidable.

We caution the reader that in the infinite time context, a decidable  model might not be computable (see Corollary
\ref{Corollary.DecidableDiagramNotComputableModel}). This is a consequence of the phenomenon in infinite time computability  that a function can have
a decidable graph without being a computable function. The classical algorithm to  compute a function from its graph relies on having an effective
enumeration of the possible values of the function, but in the infinite time context we  have no effective method to enumerate $\R$. For a purely
relational model, with no function or constant symbols in the language, however, this phenomenon  is avoided and the model is computable if and only
if its atomic diagram is decidable.

Another departure from the classical theory is that every computable  model $\mathcal A$ with underlying set contained in $\N$ is decidable.  The
point is that the infinite time algorithm can systematically check the truth  of any first order statement $\varphi$ in $\mathcal A$, given the
G\"odel code $\gcode{\varphi}$, by inductively applying the Tarski definition of  truth. If $\varphi$ has the form $\exists x\,\psi(x)$, then the
algorithm simply checks the truth of all $\psi(n)$ for $n\in A$. More generally, if an  infinite time Turing machine has the capacity for a complete
search  through the domain of a structure---for example if the domain consisted of a  writable set of writable reals---then we will be able
effectively to  carry out the Tarski definition of truth. So one might want to regard such a situation as a special or trivial case in infinite time
computable model theory. We refer to such a structure as a \emph{writable structure}; a formal definition appears on page \pageref{defn.writable}.

A theory (meaning any set of sentences in a fixed language) is {\df  computably axiomatizable} if there is a theory $T_0$, having the same
consequences as $T$, such that the set of G\"odel codes $\set{\gcode{\varphi}\st\varphi\in T_0}$ is decidable. A theory $T$ is  {\df decidable} if
the set of G\"odel codes of its consequences $\set{\gcode{\varphi}\st T\proves\varphi}$ is decidable. If the underlying language is coded in $\N$,
then every computably axiomatizable theory is decidable, because an infinite time algorithm  is easily able to search through all proofs. More
generally, if an  algorithm can write a real listing all the G\"odel codes of symbols in the language,  then it can systematically generate the
G\"odel codes of all sentences  in that language, determine which are axioms in $T_0$, and then generate a list  of all possible proofs. This shows
that any theory with a writable set  of axioms has a writable set of theorems.

\SubSection Coding with reals

We would like to view our algorithms as engaging with arbitrary countable objects, such as countable ordinals or theories, even though  formally the
machines treat only infinite binary sequences. So let us introduce a  method of coding. We regard any real $x\in\R$ as coding a relation
$\trianglelt$ on $\N$ by $i\trianglelt j$ if and only if the $\<i,j>^\th$ bit of $x$  is $1$, using a bijective pairing function $\<{\cdot},{\cdot}>$
on $\N$. For every countable ordinal $\alpha$, there is such a relation  $\trianglelt$ on $\N$ with $\<\alpha,<>\iso\<A,\trianglelt>$, where $A$  is
the field of $\trianglelt$. The set $\WO$ consists of the reals $x$ coding such well  ordered relations $\trianglelt$, and we refer to these as the
reals  coding ordinals. This is well known to be a complete $\Pi^1_1$ set of reals.  One of the early results of
\cite{HamkinsLewis2000:InfiniteTimeTM}  showing the power of infinite time Turing machines is that this set is decidable.  We sketch the proof
because the method will be useful for other purposes here.

\Theorem.({\cite[Theorem 2.2]{HamkinsLewis2000:InfiniteTimeTM}}) $\WO$  is infinite time decidable.

\Proof. Given a real $x$, we first check whether $x$ codes a linear  order $\trianglelt$, by systematically checking all instances of  transitivity,
reflexivity and anti-symmetry, in $\omega$ many steps of computation.  Assuming $\trianglelt$ is a linear order, we next attempt to find the
$\trianglelt$-least element in the field of the relation. This can be  done by placing a current guess for the least element on the scratch  tape,
and searching for a $\trianglelt$-smaller element. When such a better guess  is found, the algorithm over-writes it on the scratch tape, and also
flashes a special flag on and then off. At the next limit stage, if the flag is  on, then the guess was changed infinitely many times, and so the
real is rejected, because it does not code a well order. If the flag is off at  a limit, then the guesses stabilized on the current
$\trianglelt$-least  element, which now appears on the scratch tape. Next, the algorithm erases all mention of this element from the field of the
relation coded on the  input tape, and then continues to find (and subsequently erase) the next least  element, and so on. The algorithm should
detect limits of limit stages,  so that the scratch tape and the flag can be accordingly reset. Eventually, the  well ordered initial segment of
$\trianglelt$ is erased from the field  of the relation coded on the input tape. By detecting when the tape is empty,  the algorithm can know whether
the original real coded a well order. If  not, the algorithm will detect the ill-founded part of it, and reject at that  stage.\QED

Since $\WO$ is a complete $\Pi^1_1$ set, any $\Pi^1_1$ question reduces  to a question about $\WO$, and so we obtain:

\Corollary. Any $\Pi^1_1$ set is infinite time decidable. Hence, any  $\Sigma^1_1$ set is also decidable.\label{Corollary.Pi^1_1Decidable}

Any real $x$ can be viewed as the code of an $\omega$-sequence of reals  $\<(x)_n\st n<\omega>$ by $(x)_n(m)=x(\<n,m>)$. Thus, if we are also  given
a real $z$ coding a relation $\trianglelt$ on $\N$ of order type  $\alpha$, then any $\beta<\alpha$ is represented by some $n$ with  respect to
$\trianglelt$, and we may view $x$ as coding via $z$ an  $\alpha$-sequence $\<x_\beta\st\beta<\alpha>$ of reals. The real  $x_\beta$ is $(x)_n$,
where $n$ is the $\beta^\th$ element with respect to $\trianglelt$.

More generally, any hereditarily countable set $a$ can be coded with a  real as follows. Suppose $b$ is any countable transitive set containing  $a$
as an element, such as the transitive closure $\TC(\singleton{a})$, and  let $E$ be a relation on a subset $A\of\N$ such that there is an isomorphism
$\pi:\<A,E>\iso\<b,{\in}>$. Since this isomorphism $\pi$ must be the  Mostowski collapse of $E$, the set $a$ is determined by $E$ and the natural
number $n$ such that $\pi(n)=a$. We view the pair $\<n,E>$, coded by a real,  as representing the set $a$. Of course, a set $a$ generally has many
different codes. In analogy with $\WO$, let us define $\HC$ to be the set of such  reals coding hereditarily countable sets in this way. Given two
such codes $x$ and $y$, define $x\equiv y$ if $x$ and $y$ are codes for the same set,  and $x\in^* y$ if the set coded by $x$ is an element of the
set coded  by $y$.

\Theorem. The structure $\<\HC,{\in^*},\equiv>$ is infinite time computable.\label{Theorem.HCIsComputable}

\Proof. The elements of $\HC$ are precisely the reals coding pairs $\<n,E>$ where $E$ is a well-founded relation on some $A\of\N$, where $A$ is the
field of $E$, the natural number $n$ is in $A$, and the structure $\<A,E>$ satisfies extensionality. Thus, the set $\HC$ is $\Pi^1_1$ definable, and
hence decidable. The relation $x\equiv y$ is satisfied, where $x=\<n,E>$ and $y=\<n',E'>$ if and only  if there is an isomorphism from the part of
the field of $E$ below $n$  to the field of $E'$ below $n'$. This is a $\Sigma^1_1$ property in the codes,  and hence decidable. Similarly, the
relation $x\in^* y$ simply asserts that there is some $m$ in the field of $E'$ such that $\<n,E>\equiv\<m,E'>$,  which is also $\Sigma^1_1$, and
hence decidable.\QED

The quotient structure $\HC/{\equiv}$, under the induced relation $\in^*$, is of course isomorphic to the transitive collection $H_{\omega_1}$ of
hereditarily countable sets.

\Theorem. The satisfaction relation for hereditarily countable sets $\<b,{\in}>\satisfies\varphi[\vec a]$ is infinite time decidable, given any code
$\<n,E>\in\HC$ for b, the G\"odel code $\gcode{\varphi}$ and the code $\vec n$ of $\vec a$ with respect to
$E$.\label{Theorem.SatisfactionIsDecidable}

\Proof. This is simply an instance of the earlier remark we made, that  when an algorithm has access to the entire domain of a structure, it  can
carry out the Tarski definition of truth. In this case, the code for $b$  effectively provides the structure $\<b,{\in}>$ as a subset of $\N$.
Alternatively, one could simply observe that the satisfaction relation has complexity $\Delta^1_1$, and is therefore decidable.\QED

The constructible hierarchy of G\"odel is the transfinite hierarchy of  sets $L_\alpha$, defined by: $L_0=\emptyset$; $L_{\alpha+1}$ is the
collection of definable subsets of $L_\alpha$; for limit ordinals, $L_\eta=\Union_{\alpha<\eta}L_\alpha$. The constructible universe $L$  is the
proper class $\Union_\alpha L_\alpha$, and G\"odel proved that $\<L,{\in}>$ is a  (class) model of $\ZFC+\GCH$ and much more.

\Theorem. \
\begin{enumerate}
 \item There is an infinite time algorithm such that, on input a code for $L_\alpha$ for some countable ordinal $\alpha$, writes a code of $L_{\alpha+1}$.
 \item There is an infinite time algorithm such that, on input a code of a countable ordinal $\alpha$, writes a code of $L_\alpha$.
\end{enumerate}\label{Theorem.CodingLalpha}

\Proof. Given a code of $L_\alpha$, one systematically considers each  definition and each parameter, and by repeated applications of Theorem
\ref{Theorem.SatisfactionIsDecidable}, one can write down codes for  each of the definable subsets. This produces a code for $L_{\alpha+1}$.  Given a
code for $\alpha$, one views $\N$ as an $\alpha$-sequence of copies of  $\N$. On each copy of $\N$, the algorithm may iteratively apply the  previous
method to produce codes for the successive new elements of $L_\beta$  for each $\beta\leq\alpha$.\QED

The next theorem asserts that there is a real $c$ such that an infinite  time Turing machine can recognize whether a given real is $c$ or not,  but
no algorithm can produce $c$ on its own. This is like a person who is able  to recognize a particular song, a lost melody, when someone else sings
it, but who is unable to sing it on his or her own. The idea of the proof leads  to the concept of $L$-codes for sets and ordinals, of which we will
make extensive use later.

\Theorem Lost Melody Theorem.(\cite{HamkinsLewis2000:InfiniteTimeTM})  There is a real $c$ such that $\singleton{c}$ is infinite time  decidable, but
$c$ is not writable.\label{Theorem.LostMelodyTheorem}

\Proof. We sketch the proof from \cite{HamkinsLewis2000:InfiniteTimeTM}. Results there show that  every infinite time Turing machine computation
either halts or repeats by some countable stage. Let $\beta$ be the supremum of the stages by which all  computations of the form $\varphi_p(0)$ have
either halted or repeated. (Welch proved in \cite{Welch2000:LengthsOfITTM} that $\beta=\Sigma$, the  supremum of the accidentally writable ordinals.)
The structure $L_\beta$ is able to carry out all the computations $\varphi_p(0)$ for any length up to  $\beta$, and so the defining property of
$\beta$ is expressible in $L_\beta$. One can use the defining property of $\beta$ to show that there is a map from $\omega$ unbounded in $\beta$ that
is a definable subset of $L_\beta$. This map is therefore an element of $L_{\beta+1}$, and consequently $\beta$ is countable in $L_{\beta+1}$. So
there is some $L$-least real $c\in L_{\beta+1}$ coding a relation of order type $\beta$. This is the real we seek.

Notice that $\singleton{c}$ is decidable, because if we are given any  candidate real $c'$, we can check that it codes an ordinal $\beta'$,  and if
so, we can write down a code for $L_{\beta'+1}$, and check whether  $L_{\beta'+1}$ satisfies that $\beta'$ is the supremum of the repeat  points for
all computations $\varphi_p(0)$. This will be true if and only if  $\beta'=\beta$. Next, we check that $c'$ is the least real in
$L_{\beta'+1}=L_{\beta+1}$ coding $\beta'=\beta$. This will be true if and only if $c'=c$. So we  can decide whether any given real is $c$ or not.

Finally, $c$ is not writable, because $\beta$ is necessarily larger  than every clockable ordinal, and hence larger than every writable  ordinal. So
$\beta$ is not coded by any writable real.\QED

\Corollary. There is a function $f$ that is not infinite time computable, but whose graph is infinite time
decidable.\label{Corollary.NotComputableGraphDecidable}

\Proof. Let $f(x)=c$ be the constant function with value $c$, the lost  melody real. Since $\singleton{c}$ is decidable, we can decide the  graph of
$f$, which consists of all pairs $(x,y)$ for which $y=c$. But $f$ is not  computable, since $c\neq\varphi_p(0)$ for every program $p$.\QED

\Corollary. There is an infinite time decidable model that is not infinite time computable.\label{Corollary.DecidableDiagramNotComputableModel}

\Proof. Let ${\mathcal A}=\<\R,f>$, where $f(x)=c$ is the constant function with value $c$, given by the Lost Melody Theorem, and $\gcode{f\,}\in\N$.
This is not a computable model, because the function $f$ is not computable. Nevertheless, we  will show that the elementary diagram of $\mathcal A$
is decidable. First, we consider the atomic diagram. We can use $f(f(x))=f(x)$ to reduce the complexity of terms, and then observe that $f(x)=f(y)$
is always true and $f(x)=y$ amounts to $y=c$, which is decidable. So any atomic assertion is  decidable. To decide the full elementary diagram, we
observe that it admits the effective elimination of quantifiers down to Boolean combinations of  assertions of the form $x=c$ and $x=y$ (plus true
and false). The quantifier case essentially amounts to observing that $\exists x\,(x=c\And x\neq y)$ is  equivalent to $y\neq c$ and $\exists
x\,(x\neq c\And x\neq y)$ is  simply true. So $\mathcal A$ is decidable, but not computable, concluding the proof.

This Corollary can also be proved by using a language with a single constant symbol 0, with $\gcode{0}\in\N$. The structure ${\mathcal B}=\<\R,c>$,
interpreting $0$ as the Lost Melody real $c$, is not a computable model because the value of the constant is not computable from its G\"odel code.
But the structure $\mathcal B$ is simply an infinite model with a distinguished constant, which admits the elimination of quantifiers, and since one
can decide all statements of the form $x=c$, it follows that $\mathcal B$ has a decidable theory.\QED

The idea of the Lost Melody Theorem  provides a method of coding countable ordinals in $L$ with unique codes. Specifically, for any
$\alpha<\omega_1^L$, let $\beta$ be least above $\alpha$ such that  $\beta$ is countable in $L_{\beta+1}$, and let $c$ be the $L$-least  real in
$L_{\beta+1}$ coding a relation $\trianglelt$ on $\N$ with order type  $\beta$. The ordinal $\alpha$ is represented by some natural number $n$  with
respect to $\trianglelt$, and so we will define $\<n,c>$ to be the {\df $L$-code} of $\alpha$. Note that every ordinal $\alpha$ that is countable in
$L$ has exactly one $L$-code, since $\alpha$ determines $\beta$, which  determines $c$, which determines $\trianglelt$, which determines $n$.  Since
the $L$-code of $\alpha$ is also a code of $\alpha$ in the sense of $\HC$, we can computably determine by Theorem \ref{Theorem.HCIsComputable}
whether $\alpha<\beta$, given $L$-codes for $\alpha$ and $\beta$. And just as  with $\HC$ in this case, we can computably construct the isomorphism
from the field of the relation coding $\alpha$ to the appropriate initial  segment of the field of the relation coding $\beta$, and find the
particular natural number representing $\alpha$ with respect to the code for $\beta$.

\Lemma. The set of $L$-codes for ordinals is infinite time decidable.\label{Lemma.LcodesForOrdinalsDecidable}

\Proof. Given a real coding a pair $\<n,c>$, we can determine whether  $c$ is the code of a relation $\trianglelt$ on $\N$ that is a well  order of
some order type $\beta$. If so, we can construct a code for $L_{\beta+1}$  and check that $L_{\beta+1}$ satisfies that $\beta$ is countable and  that
the $L$-least real coding a relation of order type $\beta$ is $c$. Finally,  we can check that $L_{\beta+1}$ thinks that $\beta$ is least such that
it satisfies that $\alpha$, the ordinal coded by $n$ with respect to $\trianglelt$, is countable. If all these tests are passed, then the  pair
$\<n,c>$ is the $L$-code of $\alpha$.\QED

More generally, we have $L$-codes for any set that is hereditarily countable in $L$. Specifically, suppose that $a$ is any set that is  hereditarily
countable in $L$. Let $\beta$ be least such that $a\in L_\beta$ and  $\beta$ is countable in $L_{\beta+1}$. It follows that $L_\beta$ is  countable
in $L_{\beta+1}$, so there is some $L$-least real $c$ coding a relation  $E$ such that $\<\N,E>\iso\<L_\beta,\in>$. The set $a$ is represented by
some natural number $n$ with respect to $E$, and the {\df $L$-code} of $a$  is the pair $\<n,c>$. Let $\LC$ be the set of such $L$-codes for
hereditarily countable sets in $L$. Since these are also codes for sets in the sense  of $\HC$, it follows by Theorem \ref{Theorem.HCIsComputable}
that we may computably decide the relation $\in^*$ on the codes induced by the  $\in$ relation on the sets coded.

\Theorem. The structure $\<L_{\omega_1^L},\in>$ has an infinite time  computable presentation as
$\<\LC,\in^*>$.\label{Theorem.Lomega1HasComputablePresentation}

\Proof. The set $L_{\omega_1^L}$ is precisely $\HC^L$, the sets that  are hereditarily countable in $L$, and this is isomorphic to $\<\LC,\in^*>$ via
the $L$-codes. The $\in^*$ relation is computable on the $L$-codes,  just as in Theorem \ref{Theorem.HCIsComputable}. And the set of  $L$-codes $\LC$
is decidable just as in Lemma \ref{Lemma.LcodesForOrdinalsDecidable}.\QED

Similarly, using the $L$-codes for ordinals, we see that the structure  $\<\omega_1^L,<>$ has an infinite time computable presentation.

\Section Arithmetic on the Real Line

As a straightforward example of an infinite time computable structure, we consider the most prominent uncountable structure in mathematics, the real
line under arithmetic.

\Lemma. \label{lemma:realline} The standard structure $\RL$ of the real line under addition, multiplication, subtraction, division, and the order
relation $<$ is infinite time computably presentable.

We use ``the real line'' to describe this structure, and refer to its elements as ``points,'' because elsewhere in this paper we use the term ``real
number'' to refer to elements of $2^\omega$. Also, since division by zero is usually undefined, let us regard it as a function on the computable
domain $\RL\times (\RL-\{0\})$.

\Proof. It is straightforward to identify points $x$ on the real line uniquely with binary sequences $C\in 2^\omega$ such that $C(2n)=0$ for
infinitely many $n$ and $C(2n+1)=0$ for all but finitely many $n$ and $C\neq \<1000\cdots>$. The element $C$ corresponds to the real point
$$ (-1)^{C(0)}\left(\sum_{n=0}^\infty 2^n\cdot C(2n+1) +
\sum_{n=1}^\infty \frac{C(2n)}{2^n}\right),$$ and $C$ is called the \emph{presentation} of the real point $x$. (The condition $C\neq \<1000\cdots>$
rules out the second presentation of the point $0$ as $-0$.) The domain of our structure $\RL$ is the set of all presentations of real points, and is
decidable in infinite time, since each of the conditions can be checked in $\omega$ many steps. Also, it is easy to give a process for deciding (in
infinite time) whether two given domain elements are equal, and if not, which is larger under $<$.

Of course, all of the usual arithmetic operations on these representations, such as sum, difference, product and quotient, have complexity (much less
than) $\Delta^1_1$ in the input, and therefore, by Corollary \ref{Corollary.Pi^1_1Decidable}, these are all infinite time computable operations.
Nevertheless, for illustration let us show in moderate detail how to compute the sum $C''$ of two presentations $C$ and $C'$ of positive real points.
First, in $\omega$ many steps, we find the greatest $k>0$ such that $C(2k)=C'(2k)$, or else establish that there are infinitely many such $k$. Then
we have two cases. If there is a greatest $k$, then beyond the $k$-th bit $C$ and $C'$ complement each other perfectly, and there are only finitely
many bits left to add.
Otherwise, there are infinitely many $k$ with $C(2k)=C'(2k)$, and we build the sum from the inside, by always searching for the next greater bit $k$
with $C(2k)=C'(2k)$ and computing $C''$ up to that bit. (The point is that when $C(2k)=C'(2k)$, we know right away whether we need to ``carry'' a $1$
from $C''(2k)$ when calculating $C''(2k-2)$, even without knowing $C''(2k)$ itself yet.) This defines the entire sequence $C''$. Note that if the
representation of the sum happens to have $C''(2k)=1$ for a tail segment, then one must switch to the preferred representation by changing these bits
to $0$ and performing an additional carry.

Notice that each of the two cases could be carried out in finite time computability, producing each bit $C''(n)$ in finitely many steps, assuming
that one was given oracles presenting $C$ and $C'$. Infinite time is required only to decide which of the two cases to use, and (in the first case)
to find the greatest $k$.

Addition of two negative real points can be defined using the above algorithm conjugated by the negation map $x\mapsto -x$, which is immediately seen
to be computable. To get addition of a positive to a negative, we define subtraction of a positive real $C'$ from another one $C>C'$, by taking
finite approximations of the difference, adding them to $C'$, and checking whether each finite approximation yields a sum $>C$ or $\leq C$.

It is tempting to bypass the discussion for subtraction by saying that the difference $C-C'$ should be that domain element $D$ such that $C'+D=C$,
since we have already given a method of computing the sum of positive domain elements. However, this does not suffice to prove computability, and
indeed it illustrates a fundamental difference between the contexts of finite and infinite time: in infinite time computability, we may no longer
have such effective search procedures. Without an infinite-time-computable enumeration of the domain of $\RL$, there is no guarantee that we would
ever find the element $D$ described above, even though it must lie somewhere in the domain of $\RL$. Therefore, it is necessary to compute $D$
directly in infinite time, rather than searching for a $D$ which satisfies $C'+D=C$.  Decidability of subtraction as a ternary relation (that is,
decidability of the statement $C-C'=D$) does follow from decidability of the addition relation, which follows from computability of addition as a
function, but computability of subtraction is stronger.

For multiplication of positive domain elements $C$ and $C'$, we simply multiply $C$ by each individual bit of $C'$ (for instance, if $C'(2n)=1$, then
the product of $C$ with that bit maps each bit of $C$ $n$ places to the right) and add the results together, one by one, in $\omega^2$ many steps.
Clearly each bit on the output tape does converge to a limit, since $C(2n+1)=0$ for all but finitely many $n$, and the final output is the product of
$C$ and $C'$. This extends easily to the case of non-positive domain elements, so multiplication is computable. Finally, for division, we can check
whether the divisor is the real point $0$, and if not, we define it using the multiplication function, just as subtraction was defined using
addition. Thus division is indeed a computable function on the domain $\RL\times (\RL-\{ 0\})$. \QED

One can expand the real field $\RL$ to include all the usual functions of analysis: $e^x$, $\sqrt{x}$, $\ln x$, $\sin x$ and so on. Since (the bit
values of) these functions have complexity below $\Delta^1_1$, they are all infinite time computable by Corollary \ref{Corollary.Pi^1_1Decidable}.

Let us turn now to the subfield $\Rw$, consisting of those real points having a writable presentation. It is clear from the algorithms given in that
proof that $\Rw$ is a substructure of the real line $\RL$. Moreover, we have the following lemma.

\Lemma. \label{lemma:writables} In the infinite time context, the ordered field $\Rw$ is computably  presentable, and more generally, the ordered
field $\Rw^X$ of those real points which have presentations writable using any oracle $X\of\R$ is $X$-computably presentable. In each presentation,
there is a computable (resp.~$X$-computable) function from domain elements to the binary expansions of the real points they represent.

\Proof. The main difficulty is in getting the domain of our presentation of $\Rw$ to be decidable.  For our domain $S$, we take the set of pairs
$\langle e,c\rangle\in\omega\times\{ c\}$, where $c$ is the Lost Melody real of Theorem \ref{Theorem.LostMelodyTheorem} and the $e^\th$ infinite time
program outputs a presentation of a real point and no $e'<e$ is the index of a program outputting the same point. This is indeed a decidable domain:
given any pair, we first check whether the second element is $c$ (since the set $\{ c\}$ is decidable), and, if so, use $c$ to check the remaining
conditions, which we can now do because $c$ codes an ordinal $\alpha$ so large that every program which halts at all must halt by stage $\alpha$, as
seen in the proof of Theorem \ref{Theorem.LostMelodyTheorem}.


Given any two elements $\langle e,c\rangle$ and $\langle e',c\rangle$ of $S$, we need to compute their sum, product, difference, and quotient, and
also to compute the relation $<$. For each of the four operations, the proof of Lemma \ref{lemma:realline} gives a program $P_{e_0}$ which writes a
presentation of the resulting real point, with $e_0$ being infinite time computable uniformly in $e$ and $e'$.
So the result of the operation is the element $\langle e_1,c\rangle$, where $e_1$ is the least index of a program which outputs a presentation of the
same real point as $e_0$.  We were given $c$ itself, of course, as part of the points $\langle e,c\rangle$ and $\langle e',c\rangle$, and with $c$ it
is simple to find the least such $e_1$. Thus each operation is computable on the domain $S$. The final claim is clear, since an element of $S$
contains an algorithm for writing out a presentation of the corresponding real point, which in turn quickly yields every digit of the binary
expansion of that point. From this, the relation $<$ on $S$ is easily computed.

For $\Rw^X$, one simply relativizes the entire proof (including the choice of the Lost Melody real) to the oracle $X$. \QED

Lemma \ref{lemma:writables} shows how infinite time computable model theory differs from its finite time analogue. While we have proved that the
ordered field $\Rw$ of infinite time computable reals (i.e. the writable reals) has an infinite time computable presentation, the corresponding fact
in finite time is not true, for the finite time computable reals have no finite time computable presentation.

\begin{proposition}
\label{prop:computablereals} Let $\Rw$ be the ordered field of infinite time computable real points,
and let $\Rc$ be the ordered subfield of finite time computable real points. Then neither $\Rw$ nor $\Rc$ is finite time computably presentable (in
domain $\N$), but both are computably presentable in infinite time.
\end{proposition}

Since the rational ordered field $\Q$ embeds uniquely and densely into $\RL$, it follows that every ordered subfield $\F$ of $\RL$, such as $\Rw$ or
$\Rc$, embeds uniquely into our presentation of $\RL$. We show next that this unique embedding is computable.

\Lemma. If\/ $\F$ is any computable presentation of an ordered subfield of $\RL$, then the unique embedding of $\F$ into $\RL$ is computable.

\Proof. Given any $x\in\F$, we may use the computable functions of $\F$ to systematically compute the $\F$-representations of the rationals
$\frac{m}{2^n}$, and make comparisons of these rationals with $x$ using the order of $\F$. This allows us to know the binary representation of $x$,
and therefore also the representation of $x$ in our presentation of $\RL$. Thus, we have computed the unique embedding of $\F$ into our presentation
of $\RL$.\QED

\Proof Proof of Proposition \ref{prop:computablereals}. An infinite time computable presentation of $\Rw$ was shown above to exist, and $\Rc$ is an
infinite time computable subset of the domain, since infinite time Turing machines can easily simulate finite time ones.

If $\F$ were a finite time computable presentation of $\Rc$, then given any element $x\in\F$, we could compute the $n$-th digit of the binary
expansion of the real point corresponding to $\F$, in finite time and uniformly in $x$ and $n$. If $\F \cong \Rc$, this would give a simultaneous
uniform finite time computation of all finite time computable sets, which of course is impossible. If $\F \cong \Rw$, then it would give a
simultaneous uniform finite time computation of all infinite time writable reals, which again is easy to diagonalize against. This completes the
proof of Proposition \ref{prop:computablereals}. \QED

We note that the same diagonalization against finite time computable presentations of all finite time computable sets can be used to show that there
is no infinite time writable presentation of all infinite time writable reals. Therefore we ask how it is that $\Rw$ is infinite time computably
presentable. The answer is that while the domain of the presentation of $\Rw$ is a countable decidable set, it is not the image of $\omega$ under any
infinite time computable function.  The use of the lost-melody real $c$ in the domain of $\Rw$ makes this clear, and indeed, without using $c$ or a
similar element, we could not decide in infinite time which programs output
infinite time computable reals.  

A concise statement of the foregoing argument is to say that there is no writable presentation of $\Rw$, even though there is a computable
presentation. \label{defn.writable} A \emph{writable structure} is an infinite time computable structure $\mathcal{A}$ such that there exists a
single writable real $r\in 2^\omega$ whose first row $r^{[0]}$ codes the entire atomic diagram of $\mathcal{A}$ and whose remaining  rows name all
elements of the domain of $\mathcal{A}$.  That is,
$$ r^{[0]} = \{\gcode{\varphi}~:~\varphi\in D_a(\mathcal{A})\}$$
$$ \text{dom}(\mathcal{A}) = \{r^{[n]}~:~n\in\omega -\{ 0\}\}.$$
(An equivalent definition requires that $r^{[n]}\neq r^{[m]}$ whenever $0<n<m<1+|\mathcal{A}|$, and $r^{[m]}=0$ if $m > |\mathcal{A}|$.) We assume
for these purposes that the language of $\mathcal{A}_A$ is also coded into $\omega$, with $2n-1$ coding the constant symbol for the element named as
$r^{[n]}$. Thus we have a computable enumeration of the elements of $\mathcal{A}$, from which it is immediate that the complete diagram of
$\mathcal{A}$ is infinite time decidable. Since they allow computable searches of the entire domain, writable structures behave something like an
analogue to the finite structures in the classical theory.

Let us conclude this section with a brief generalization.
Let $\Ra$ be the structure of the real points having an accidentally writable presentation, and similarly, let $\Rev$ consist of those having an
eventually writable presentation.

\Theorem. $\Rf\elesub\Rw\elesub\Rev\elesub\Ra\elesub \RL$.

\Proof. The point is that each of these structures is a real closed ordered field. Before explaining this, let us first iron out a wrinkle with
$\Ra$. In order to see even that this structure is closed under addition, it is useful to know that the set of accidentally writable reals is closed
under pairing. To see this, consider the algorithm that simulates all programs on input $0$, and for each accidentally writable real $x$ observed
during this master simulation, the algorithm starts another master simulation that produces all accidentally writable reals $y$ that appear before
the first appearance of $x$. Then, for each such $y$, our main algorithm writes a real coding the pair $\<x,y>$ on the scratch tape. This algorithm
shows that if $x$ and $y$ are accidentally writable, then the pair $\<x,y>$ is also accidentally writable. Using this and the observations of Lemma
\ref{lemma:realline}, it now follows that $\Ra$ is a field.

Each of the fields is closed under square roots for its positive elements, since the digits of the square root can be systematically computed. Also,
for any odd degree polynomial, one can use successive approximations (for example, by Newton's method) to find a computable root. Since the theory of
real closed fields is model complete, the theorem now follows.\QED

One can naturally extend this theorem by oracles and have a rich lattice of relatively computable subfields of $\RL$. Each of the extensions in the
theorem is strict, by \cite[Theorem 6.15]{HamkinsLewis2000:InfiniteTimeTM}, and it follows that each is a transcendental extension of the previous.
Finally, we observe that $\Rw$ can have no writable transcendence basis over $\Q$ or $\Rf$, since then we would be able to produce a writable list of
all writable reals, which we have observed is impossible by a simple diagonalization. Similarly, $\Rev$ has no eventually writable transcendence
basis over $\Rw$ and $\Ra$ has no accidentally writable transcendence basis over $\Rev$.

\Section The Infinite time computable completeness theorem

The Completeness theorem asserts that every consistent theory has a  model. The finite time effective version of this asserts that any  finite time
decidable theory has a finite time decidable model. And in the infinite  time context, at least for languages coded in $\N$, this proof goes  through
without any hitch. In fact, the infinitary context gives a slightly  stronger result:

\Theorem. In the infinite time context, if\/ $T$ is a consistent theory  in a computable language coded in $\N$ and $T$ has a computable
axiomatization, then $T$ has a decidable computable model. In fact, such a theory has a  model coded by a writable
real.\label{Theorem.CompletenessTheoremLanguageInN}

\Proof. The point is that the classical Henkin construction is effective for infinite time Turing machines. Note that if $T$ has a  computable
axiomatization in a language coded in $\N$, then it is actually  decidable, since the infinite time Turing machines can search through  all proofs in
$\omega$ steps. We may assume that there is an infinite supply of new  constant symbols, by temporarily rearranging the G\"odel codes of the  symbols
in the original language if necessary. Enumerate the sentences in the expanded langauge as $\<\sigma_n\st n\in\N>$, and build a complete consistent
Henkin theory in the usual manner: at stage $n$, we add $\sigma_n$, if this is  consistent with what we have already added to $T$, or else
$\neg\sigma_n$, if it is not. Since $T$ is decidable, this is computable. In addition, if  $\sigma_n$ has the form $\exists x\,\varphi(x)$ and we
added it to the  theory, then we also add $\varphi(c)$ for the first new constant symbol $c$ that has not yet been considered. The result of this
construction is a  complete consistent Henkin theory $\bar T$ extending $T$. The theory $\bar T$ is  decidable, because for any $\sigma$, the
infinite time algorithm can  run the construction until $\sigma$ is considered, and answer accordingly as it  was added to $\bar T$ or not. As usual,
we may use the Henkin constants  to build a model of $T$. Specifically, let $c\equiv d$ if $\bar T\proves c=d$, and define the $R([\vec c])\iff\bar
T\proves R(\vec c)$ and $f([\vec c])=[d]\iff\bar T\proves f(\vec c)=d$. The classical induction shows  that the resulting structure $M_{\bar T}$ of
equivalence classes  satisfies $\varphi([\vec c])$ if and only if $\bar T\proves\varphi(\vec c)$, so  this is a model of $T$. Finally, for any
constant symbol $d$, one may  compute the (numerically) least element of $[d]$ by simply testing each of the smaller constants $c$ to determine
whether $c\equiv d$. Thus, by replacing each equivalence class with its least member, we construct a computable  presentation of $M_{\bar T}$. Since
the underlying set of this model is  contained in $\N$, an algorithm can write down the entire structure as a writable  real.\QED

Many theories, including some very powerful theories, have infinite  time computable axiomatizations, and so this result provides numerous
interesting decidable models. For example, the theory of {\df true arithmetic}  $\TA=\Th(\<\N,+,\cdot,0,1,<>)$ is infinite time decidable, because
arithmetic truth is infinite time decidable, and so the theory $\TA+\set{n<c\st n\in\N}$  is a computable axiomatization of the theory of the
nonstandard models of true arithmetic. Similar observations establish:

\Corollary. There are infinite time decidable computable nonstandard  models of the theories \PA, \TA, \ZFC, $\ZFC\ +$ large cardinals, and  so on,
provided that these theories are consistent.

The infinite time realm, therefore, lies considerably beyond the computable models of the finite time theory. What is more, as we have  emphasized,
the infinite time context allows for uncountable computable models and  uncountable languages, which cannot be coded in $\N$. So Theorem
\ref{Theorem.CompletenessTheoremLanguageInN} doesn't tell the full story. In the general context, where languages are coded in the reals,  we ask
whether the full infinite time analogue of the Completeness Theorem  holds:

\Question. Does every consistent infinite time decidable theory have an  infinite time decidable model? Does every such theory have an infinite  time
computable model?\label{Question.ComputableCompleteness}

One of the convenient features of the classical theory, when working  with a language coded in $\N$, is that one can enumerate the function, relation
and constant symbols of the language $s_0$, $s_1$, $\ldots$ in such a  way that from any symbol $s_n$, one can reconstruct the list $\<s_m\st  m\leq
n>$ of prior symbols. This is a triviality in the context of computable  languages coded in $\N$, because we simply enumerate the symbols in the
order of their G\"odel codes. Given any such code, one simply tests all the smaller natural numbers in turn to discover the list of prior codes for
symbols. But in the uncountable context, a computable representation of a language  may not have this feature. Let us therefore define that a
computable representation $\mathcal L$ of a language is computably {\df well  presented} if there is an enumeration $\<s_\alpha\st\alpha<\delta>$ of
all of the function, relation and constant symbols of the language, for some $\delta\leq\omega_1$, such that from any $\gcode{s_\alpha}$, we can
(uniformly, in infinite time) compute a code for the sequence $\<\gcode{s_\beta}\st\beta\leq\alpha>$  of prior symbols. In this case, we can prove
the infinite time  computable analogue of the Completeness Theorem.

\Theorem. Every consistent infinite time decidable theory in a computably well presented language has an infinite time decidable model in this
language.\label{Theorem.CompletenessTheoremIfWellPresented}

We begin with a few preliminary lemmas. Let us say that a computable  presentation $\mathcal L$ of a language admits a {\df computably  stratified
enumeration} of formulas if there is an enumeration of all $\mathcal  L$-formulas $\<\varphi_\alpha\st\alpha\leq\delta>$, for some
$\delta\leq\omega_1$, such that from the G\"odel code $\gcode{\varphi_\alpha}$, one can  (uniformly in infinite time) compute a real coding the
sequence $\<\gcode{\varphi_\beta}\st\beta\leq\alpha>$ of G\"odel codes of the  prior formulas.

\SubLemma. If a language $\mathcal L$ is computably well presented,  then it admits a computably stratified enumeration of
formulas.\label{Lemma.WellPresentedIsStratified}

\Proof. Suppose that a language $\mathcal L$ is computably well presented by the  enumeration $\<s_\alpha\st\alpha\leq\delta>$. Given a well-ordered
list of function, relation and constant symbols, one can systematically produce  a list of all formulas in that language, as follows. The first
$\omega$ many formulas are those not using any of the symbols; the next $\omega$ many  formulas are those using the first symbol only; the next
$\omega$ many formulas use the second symbol and possibly the first. There is a (finite time)  computable list of countably many first order formula
templates, with  holes for the function, constant and relation symbols, and the actual formulas  are obtained by plugging codes for actual function,
relation and  constant symbols (of the appropriate arity) into those holes. From the presentation of  the symbols, we systematically generate a list
of all finite sequences of the symbols, and from these and the templates, one can generate the list of  all formulas. We therefore generate the
formulas in blocks of length  $\omega$, and all formulas in the $\alpha^\th$ block are required to use the  symbol $s_\alpha$ and may use earlier
symbols. This defines the  enumeration of the formulas $\<\varphi_\alpha\st\alpha\leq\gamma>$.

Given any formula $\gcode{\varphi}$, we can inspect it for the symbols  $s$ that appear in it, and from each $\gcode{s}$, we can generate the
corresponding list of prior symbols $\<\gcode{s_\beta}\st\beta\leq  \alpha>$, where $s=s_\alpha$. By comparing the lengths of these  sequences, we
can tell which symbol was the last to appear in the enumeration of $\mathcal  L$. For this maximal $\alpha$, we know that $\varphi$ appears in the
$\alpha^\th$ block of formulas. From the list of symbols $\<\gcode{s_\beta}\st  \beta\leq\alpha>$, we can regenerate the list of formulas up to and
including the $\alpha^\th$ block of formulas, thereby producing the prior list of  formulas $\<\gcode{\varphi_\xi}\st \xi\leq\eta>$, where
$\varphi=\varphi_\eta$.\QED

A fundamental construction of first order logic is to expand a language  by adding infinitely many new constant symbols. In the context of computable
model theory, whether finite or infinite time, if the presentation of a  language $\mathcal L$ already uses all the available G\"odel codes,  then
one is forced to consider translations of the language in order to free up  space in the G\"odel codes to represent the expanded language. For
example, even in the finite time context, if one has a model in a language with  infinitely many constant symbols, and the G\"odel codes of the
symbols  already use up all of $\N$, then in order to add constants to the language one seems  forced to use a translation of the language. A given
language can have  many different computable presentations, and in general these may not be computably equivalent. For two presentations of the
language, there may  be no computable method of translating symbols or formulas from one  representation to the other. (And this phenomenon occurs
already in the  finite time context.) In the infinite time context, where we represent symbols with  real numbers, this phenomenon can occur even in
finite languages, since the G\"odel codes for a symbol may be reals that are incomparable in the  infinite time Turing degrees. If we have two
computable presentations  ${\mathcal L}$ and ${\mathcal L}'$ of a language, and it happens that there is a  computable function mapping every
${\mathcal L}'$ code for a symbol to  the ${\mathcal L}$ code for the same symbol, then we will say that  ${\mathcal L}'$ is a {\df computable
translation} of ${\mathcal L}$. In  such a case, syntactic questions about ${\mathcal L}'$ can be computably reduced to  syntactic questions about
${\mathcal L}$. This relation is not  necessarily symmetric (because in the infinite time context, a function can be  computable without its inverse
being computable). If both languages are  computable translations of each other, we say that the languages are computably isomorphic translations.

\SubLemma. If a language $\mathcal L$ is computably well presented,  then there is a computably isomorphic translation of it to a well  presented
language ${\mathcal L}_0$, preserving the order of the enumeration of  symbols, and a well presented expansion ${\mathcal L}_1$ of ${\mathcal  L}_0$
containing $\omega$ many new constant symbols $c_s^n$ for every symbol  $s$ of $\mathcal L$, such that from $\gcode{s}$ and $n$ one can uniformly
compute $\gcode{c_s^n}$ and conversely.\label{Lemma.AddingConstants}

\Proof. For each symbol $s$ of $\mathcal L$, let its code in ${\mathcal  L}_0$ be obtained by simply adding a $0$ to the front of $\gcode{s}$ in
$\mathcal L$. For ${\mathcal L}_1$, the code of the constant symbol  $c_s^n$ is obtained by adding $n+1$ many $1$s plus $0$ to the front of
$\gcode{s}$ in $\mathcal L$. Thus, from $\gcode{s}$ in $\mathcal L$ we can easily  compute every $\gcode{c_s^n}$ and $\gcode{s}$ in ${\mathcal L}_1$
and  vice versa. So it is clear that ${\mathcal L}_0$ is a computably isomorphic  translation of the language $\mathcal L$. The enumeration of the
symbols of ${\mathcal L}_1$ simply replaces each symbol $s$ of $\mathcal L$ with the block of  symbols $s$, $c_s^0$, $c_s^1$, and so on. From any of
these symbols, we  can reconstruct the prior list of symbols in $\mathcal L$, and from those  symbols we can reconstruct the corresponding constant
symbols, so as to  generate the prior list of symbols in ${\mathcal L}_1$.\QED

\Proof Proof of Theorem \ref{Theorem.CompletenessTheoremIfWellPresented}. We carry out the proof of Theorem
\ref{Theorem.CompletenessTheoremLanguageInN} in this more general  context. Suppose that $T$ is a computably axiomatized consistent theory  in the
well presented language $\mathcal L$. Let ${\mathcal L}'$ be the well  presented language of Lemma \ref{Lemma.AddingConstants}, with  infinitely many
new constant symbols for each symbol of $\mathcal L$. Because it is well  presented, this expanded language has a computably stratified enumeration
$\<\gcode{\varphi_\alpha}\st\alpha<\delta>$ of formulas.  We assume that this language is enumerated in the manner of Lemma
\ref{Lemma.WellPresentedIsStratified}, in blocks of length $\omega$ containing all formulas with a given symbol and earlier symbols. Because we
arranged that every symbol $s$ of $\cal L$ gives rise to an infinite list of new constant symbols $c_s^n$, we may arrange that from any
$\gcode{\varphi_\alpha}$, we may uniformly compute the code of a distinct new constant symbol $c$ not appearing in any earlier $\varphi_\beta$.

We now recursively build the theory $\bar T$ in stages: at stage $\alpha$, if $\varphi_\alpha$ is a sentence, then we add it to $\bar T$ if this
remains consistent, otherwise we add $\neg\varphi_\alpha$. In addition, if $\varphi_\alpha$  is a sentence of the form $\exists x\,\psi(x)$ and we
had added it to $\bar T$, then we also add a sentence of the form $\psi(c)$, where $c$ is the distinct new constant symbol which has not yet appeared
in any earlier formula. The usual model theoretic arguments show that $\bar T$ is a complete consistent Henkin theory extending $T$.

We argue that $\bar T$ is decidable. Given any ${\mathcal L}'$ formula  $\varphi_\alpha$, we may use the computable stratification to write down a
code of $\<\gcode{\varphi_\beta}\st\beta\leq\alpha>$. From this, we may computably reconstruct $\bar T$ up to stage $\alpha$. The question of whether
to add $\varphi_\beta$ or $\neg\varphi_\beta$ at stage $\beta$ reduces to a  question about whether the theory constructed up to stage $\beta$ proves
$\neg\varphi_\beta$ or not. But since the algorithm has a real coding  the theory constructed up to stage $\beta$, it can computably enumerate  all
finite combinations of the formulas it is committed to adding to $T$, and check whether $T$ proves that any of those finite combinations of  formulas
proves $\neg\varphi_\beta$. This is a decidable question, since $T$ is  decidable and we may computably translate between the languages  $\mathcal L$
and ${\mathcal L}'$. Thus, $\bar T$ is computable.

Next, we build a decidable model of $\bar T$. Define the equivalence  relation $c\equiv d\iff\bar T\proves c=d$, and from each equivalence class
$[c]$, select the constant $c_{s_\alpha}^n$ such that the pair $\<\alpha,n>$ is lexicographically least, where $s_\alpha$ is the $\alpha^\th$ symbol
in the original presentation of $\cal L$. The set of such least constants is decidable, because from any constant $c_{s_\alpha}^n$ we may construct
the list of prior symbols, and therefore the $\omega$-blocks of the symbols in $\cal L'$, and therefore all the corresponding formulas
$\varphi_\beta$ containing only those symbols. By reconstructing the theory $\bar T$ up to that point, we can tell whether $\bar T$ proves
$c_{s_\alpha}^n=c_{s_\xi}^m$ or not, for any $\xi<\alpha$. So the set of such least representatives is decidable. We may now impose the usual
structure on these representatives, to get a decidable model of $\bar T$. Since we have a computable isomorphism of $\mathcal L'$ with $\mathcal L$,
it is no problem to translate between the two languages, and so we may use the original language presentation $\mathcal L$ when imposing this
structure, resulting in a decidable model of $T$ in the original language $\mathcal L$, as desired.\QED

\Theorem. If $V=L$, then every consistent infinite time decidable theory has an infinite time decidable model, in a computable translation of the
language.\label{Theorem.CompletessTheoremInL}

\Proof. The first step is to translate to a computably well presented  language.

\SubLemma. If $V=L$, then every computably presented language has a  computable translation to a computably well presented
language.\label{Lemma.WellPresentedTranslationsIfV=L}

\Proof. Assume $V=L$ and suppose that $\mathcal L$ is a computably presented language. Let $S\of\R$ be the corresponding computable set of  G\"odel
codes for the function, relation and constant symbols of $\mathcal L$.  Let $\<s_\alpha\st\alpha<\delta>$ be the enumeration of the elements of  $S$
in order type $\delta\leq\omega_1$, using the canonical $L$-ordering of  $\R^L$. For each $\alpha<\delta$, let $\gamma_\alpha$ be the smallest
countable ordinal above $\alpha$ such that $L_{\gamma_\alpha}$ satisfies  ``$\omega_1$ exists'' and $s_\beta$ exists for every $\beta\leq\alpha$.  By
this latter assertion, we mean that for every $\beta\leq\alpha$, the structure  $L_{\gamma_\alpha}$ computes that $S$ has at least $\beta$ many
elements in the $L$-order. Notice that because it satisfies ``$\omega_1$ exists,'' this  structure correctly computes all infinite time computations
for input  reals that it has. Therefore, it correctly computes $S\intersect  L_{\gamma_\alpha}$, which has $\<s_\beta\st\beta\leq\alpha>$ as an
initial segment in the $L$ order. In particular, $\<s_\beta\st\beta\leq\alpha>\in L_{\gamma_\alpha}$. Let $t_\alpha$ be the $L$-code of the pair
$\<\alpha,\gamma_\alpha>$. We will use $t_\alpha$ to represent the symbol coded by $s_\alpha$ in  $\mathcal L$. Denote this new translation of the
language by ${\mathcal  L}'$.

First, we observe that the set $\set{t_\alpha\st\alpha\leq\delta}$ is  decidable. Given any real $t$, we can check if it is the $L$-code of a  pair
of ordinals $\<\alpha,\gamma>$, and if so, whether $\gamma$ is least such  that $L_\gamma$ satisfies ``$\omega_1$ exists'' and $s_\beta$ exists  for
every $\beta\leq\alpha$. If so, then we accept $t$. Necessarily, in this case  $t=t_\alpha$. These questions are all decidable, because we know how
to recognize an $L$-code for a pair of ordinals, and given the code of an  ordinal $\gamma$ we can construct a code of $L_\gamma$, and then check the
truth of any statement in that structure by Theorem  \ref{Theorem.SatisfactionIsDecidable}.

What is more, from $t_\alpha$ we can construct all earlier $t_\beta$  for $\beta\leq\alpha$, because with an $L$-code for $\gamma$ we can  look for
the least $\gamma'\leq\gamma$ such that $L_{\gamma'}$ satisfies  ``$\omega_1$ exists'' and $s_\xi$ exists for all $\xi\leq\beta$. Thus,  our new
language is computably well presented via ${\mathcal L}'$. Finally, ${\mathcal L}'$  is a computable translation of $\mathcal L$ because from
$t_\alpha$ we can compute $s_\alpha$.\QED

We remark that the translation from $\mathcal L$ to ${\mathcal L}'$,  while perhaps not a computably isomorphic translation, is nevertheless
relatively mild. Specifically, from $s_\alpha$ and any code for a sufficiently  large ordinal, one can compute $t_\alpha$. In this sense, the two
representations of the language are close.

We now complete the proof of Theorem \ref{Theorem.CompletessTheoremInL}. Assume $V=L$ and suppose that $T$  is a consistent decidable theory in a
language $\mathcal L$. (By testing whether certain tautologies are well  formed, it follows that the language itself is computable.) By Lemma
\ref{Lemma.WellPresentedTranslationsIfV=L}, there is a computable translation of $\mathcal L$ to a well presented language ${\mathcal  L}'$. Let $T'$
be the corresponding translation of $T$ into this translated language.  Note that $T'$ remains decidable in ${\mathcal L}'$, because the question
$T'\proves\sigma'$ computably reduces to a question of the form  $T\proves\sigma$, which is decidable. By Theorem
\ref{Theorem.CompletenessTheoremIfWellPresented}, the theory $T'$ has a  decidable model, as desired.\QED

So it is at least consistent with \ZFC\ that the infinite time computable Completeness Theorem holds, if one allows computable translations of the
language, and in this sense one may consistently hold a positive answer  to Question \ref{Question.ComputableCompleteness}. Does this settle the
matter? No, for we will now turn to negative instances of the completeness  theorem. The fact is that in some models of set theory, there are
consistent decidable theories having no decidable model, and so the infinitary  computable Completeness Theorem is actually independent of \ZFC.

\Theorem. It is relatively consistent with \ZFC\ that there is an infinite time decidable theory, in a computably presented language,  having no
infinite time computable or decidable model in any translation of the  language (computable or not).\label{Theorem.CompletenessTheoremCanFail}

This theorem relies on the following fact from descriptive set theory. For a proof, see \cite[Theorem 25.23]{Jech:SetTheory3rdEdition}.

\SubLemma.(Mansfield-Solovay) If $A\of\R$ is $\Sigma^1_2$ and $A\not\of L$, then $A$ contains a perfect subset.\label{Lemma.MansfieldSolovay}

The crucial consequence for us will be:

\SubLemma. If $\omega_1^L$ is countable and the \CH\ fails, then there are no $\Sigma^1_2$ sets of size $\omega_1$. Hence, under these hypotheses,
there are also no decidable sets or semi-decidable sets of size $\omega_1$.\label{Lemma.NoDecidableSetsSizeOmega1}

\Proof. Every decidable or semi-decidable set $A\of\R$ is $\Delta^1_2$ and hence  $\Sigma^1_2$. If $\omega_1^L$ is countable and $A\of L$, then $A$
is countable. If $A\not\of L$, then by Lemma \ref{Lemma.MansfieldSolovay} it contains a perfect subset, and hence has cardinality $2^\omega$. Under
$\neg\CH$, this excludes the possibility that $A$ has cardinality $\omega_1$.\QED

\Proof Proof of Theorem \ref{Theorem.CompletenessTheoremCanFail}. Suppose that $\omega_1^L$ is countable and the \CH\ fails. An elementary forcing
argument shows that this  hypothesis is relatively consistent with \ZFC. Lemma \ref{Lemma.NoDecidableSetsSizeOmega1} now shows that there are no
$\Sigma^1_2$ sets of size $\omega_1$. Consider the following theory, in the language with a constant $c_x$ for every $x\in\WO$ (for simplicity, let
$\gcode{c_x}=x$), a binary relation $\equiv$ and a function symbol  $f$. The theory $T$ is the atomic diagram of the structure $\<\WO,\equiv>$, where
$\equiv$ is the relation of coding the same ordinal, together  with the axiom asserting that $f$ is a choice function on the equivalence classes.
That is, $T$ contains all the atomic facts that are true about the constants $c_x$ for $x\in\WO$, plus the assertion ``$x\equiv f(x)$ and $x\equiv
y\implies  f(x)=f(y)$.'' This theory is computably axiomatizable, because $\equiv$  is a decidable relation on $\WO$. So as a set of sentences, the
axioms of  $T$ are decidable.

But actually, the theory $T$ is fully decidable. First, we observe that  it admits elimination of quantifiers. The point is that $T$ is  essentially
similar to the theory of an equivalence relation with infinitely many  equivalence classes, all infinite. Note that the theory $T$ implies
$f(f(x))=f(x)$, and $x\equiv f(y)$ is the same as $x\equiv y$. Also,  $x=f(y)$ is equivalent to $x=f(x)\And x\equiv y$. By combining these reductions
with the usual induction, it suffices to eliminate quantifiers from  assertions of the form $\exists x\, x\equiv y\And x\not\equiv z\And x=f(x)$ and
$\exists x\, x\equiv y\And x\not\equiv z\And x\neq f(x)$. But these are  both equivalent to $y\not\equiv z$, since in the former case one may use
$x=f(y)$, and in the latter case some $x$ equivalent to $y$, other than  $f(y)$. This inductive reduction provides a computable method of finding,
for any given formula, a quantifier-free formula that is equivalent to it under $T$. The point now is that any quantifier-free sentence is a Boolean
combination of assertions about the constants $c_x$ of the form $c_x\equiv c_y$, $c_x=c_z$ and $f(c_x)=c_y$. The first two of these are computable,
since they are equivalent to $x\equiv y$ and $x=z$, respectively. The assertion  $f(c_x)=c_y$ is false if $x\not\equiv y$, which is computable, and
otherwise it is not settled by $T$, since there are models of $T$ where $f(c_x)$ is any  desired $c_y$ with $y\equiv x$. For any finite list of
constants $c_y$,  it is consistent that $f(c_x)$ is equal to any of them (at most one of them),  provided $x\equiv y$, or none of them. Because of
this, we can  computably decide whether $T$ proves any given quantifier-free assertion in the  language of $T$. So $T$ is decidable.

Finally, suppose towards contradiction that $T$ has a computable or decidable model  $M=\<A,\equiv^M,f^M,c_x^M>_{x\in\WO}$. In this case, both the
graph of $f$ and the relation $z=c_x^M$ are decidable, and so the set $\set{f(c_x^M)\st x\in \WO}$ has complexity $\Sigma^1_2$. But this set also has
cardinality $\omega_1$, contradicting Lemma \ref{Lemma.NoDecidableSetsSizeOmega1}. So $T$ can have no computable or decidable model under these set
theoretic hypotheses. Since the set theoretic hypotheses are relatively consistent with \ZFC, it is relatively consistent with \ZFC\ that there is an
infinite time decidable theory with no computable or decidable model.\QED

\Corollary. The infinite time computable Completeness Theorem is independent of \ZFC.

For this corollary we take the infinite time computable Completeness  Theorem to be the assertion: every consistent decidable theory in a  computably
presented language has a decidable model in a computable translation of  the language.

\Section The infinite time computable Lowenheim-Skolem Theorem

The classical Lowenheim-Skolem Theorem has two parts: the upward theorem asserts that every infinite model has arbitrarily large elementary
extensions, in every cardinality at least as large as the original model and the  language; the downward theorem asserts that every infinite model
has  elementary substructures of every smaller infinite cardinality at least as large as the language. Here, of course, we are interested in the
infinite time computable analogues of these assertions, which concern computable or decidable  models.

\Question. Does every infinite time decidable model have an infinite time decidable elementary extension of size continuum?
\label{Question.UpwardLowenheimSkolem}

\Question. Does every infinite time decidable infinite model (in a language coded in $\N$, say) have a countable infinite time decidable  elementary
substructure?\label{Question.DownwardLowenheimSkolem}

These questions have many close variants, depending, for example, on  whether the models are decidable or computable, and on whether the  languages
or models are well presented or not. One could ask in Question  \ref{Question.UpwardLowenheimSkolem} merely for a proper elementary  extension, or
for an uncountable extension, rather than one of size continuum, and in  Question \ref{Question.DownwardLowenheimSkolem}, merely for a proper
elementary substructure rather than a countable one (when the original model is  uncountable). We regard all such variations as infinite time
computable  analogues of the Lowenheim-Skolem Theorem.

If the Continuum Hypothesis fails badly, then it is too much to ask for  computable models of every cardinality between $\omega$ and $2^\omega$. To
be sure, this is clearly impossible if the continuum is too large (if  $2^\omega\geq\aleph_{\omega_1}$), for in this case there would be  uncountably
many such intermediate cardinalities but only countably many decidable  models. More importantly, however, Lemma \ref{Lemma.MansfieldSolovay} shows
that there can be no decidable sets of cardinality strictly between $\omega_1^L$ and $2^\omega$. Thus, the possible cardinalities of decidable sets
of reals are: finite, countable, $\omega_1^L$ and $2^\omega$.

We do not know the full answers to either of the questions above, although we do know the answers to some of the variants. For the upward  version,
if a model is well presented, then we can find an infinite time decidable  proper elementary extension (see Theorem
\ref{Theorem.UpwardLowenheimSkolemIfWellPresented}); if $V=L$, then we  can arrange this extension to be uncountable (see Theorem
\ref{Theorem.UpwardLowenheimSkolemIfV=L}). So it is consistent that the  upward Lowenheim-Skolem Theorem holds. For the downward version, if an
uncountable decidable model is well presented, then we can always find  a countable decidable elementary substructure (see Theorem
\ref{Theorem.DownwardLowenheimSkolemIfWellPresented}); but if one broadens Question \ref{Question.DownwardLowenheimSkolem} to the case of computable
models, rather than decidable models, then we have a strong negative  answer, for there is a computable structure on $\R$ having no computable proper
elementary substructures (see Theorem \ref{Theorem.ComputableModelOnRwithNoElementarySubstructures}).

In analogy with well presented languages, let us define that an infinite time computable model ${\mathcal A}=\<A,\cdots>$ is {\df well  presented} if
the language of its elementary diagram is well presented. This means  that there is an enumeration $\<s_\alpha\st\alpha<\delta>$, for some
$\delta\leq\omega_1$, including every G\"odel code for a symbol in the  language and every element of $A$, such that from $s_\alpha$ one can  compute
a code for $\<s_\beta\st\beta\leq\alpha>$. The models produced in the  computable Completeness Theorem
\ref{Theorem.CompletenessTheoremIfWellPresented}, for example, have this property.

\Theorem. If $\mathcal A$ is a well presented infinite time decidable  infinite model, then $\mathcal A$ has a proper elementary extension  with an
infinite time decidable presentation.\label{Theorem.UpwardLowenheimSkolemIfWellPresented}

\Proof. Let $T$ be the elementary diagram of $\mathcal A$, in a well  presented language. Let ${\mathcal L}'$ be the language of $T$ together  with
new constants, as in Lemma \ref{Lemma.AddingConstants}. Let $T'$ be the  theory $T$ together with the assertion that these new constants are not
equal to each other or to the original constants. Since $T$ is decidable, it is easy to see that $T'$ is decidable, since any question about whether
$T'$ proves an assertion about the new constants can be decided by replacing  them with variables and the assumption that those variables are not
equal. Thus, by Theorem \ref{Theorem.CompletenessTheoremIfWellPresented}, there is an infinite time decidable model of $T'$. Such a model provides a
decidable presentation of a proper elementary extension of $\mathcal A$.\QED

\Theorem. If\/ $V=L$, then every infinite time decidable infinite model  $\mathcal A$ elementarily embeds into an infinite time decidable model  of
size the continuum, in a computable translation of the language.\label{Theorem.UpwardLowenheimSkolemIfV=L}

\Proof. Assume $V=L$ and suppose that $\cal A$ is an infinite time decidable infinite model. We may assume, by taking a computably isomorphic copy of
the language, that all the G\"odel codes of symbols and elements in  $\cal A$ begin with the digit $0$. So there are continuum many  additional
codes, beginning with $1$, that we use as the G\"odel codes of new constant  symbols. If $T$ is the elementary diagram of $\cal A$, then let $T'$ be
$T$ together with the assertion that these new constants are not equal. The  theory $T'$ is decidable, because any question about whether $T'$ proves
an assertion reduces to a question about whether $T$ proves an assertion  about some new arbitrary but unequal elements. This can be decided by
replacing those new constant symbols with variable symbols plus the assertion  that they are distinct. Thus, by Theorem
\ref{Theorem.CompletessTheoremInL}, there is a decidable model $\cal A'\satisfies T'$. The model $\cal A'$ has  size continuum because of the
continuum many new constants we added,  and $\cal A$ embeds elementarily into $\cal A'$ because $\cal A'$ satisfies the  elementary diagram of $\cal
A$.\QED

We note that the graph of the elementary embedding of $\cal A$ into  $\cal A'$ is infinite time decidable, because from the code of a symbol  in the
expanded language, one can compute the code of the corresponding symbol in the original language. There seems little reason to expect in general that
this embedding should be a computable function, and it cannot be if the original presentation was not well  presented.

Let us turn now to the infinite time computable analogues of the downward Lowenheim-Skolem Theorem.

\Theorem. If $\mathcal A$ is an uncountable well presented infinite  time decidable model in a language coded by a writable real, then there  is an
infinite time decidable, countable elementary substructure ${\mathcal  B}\elesub{\mathcal A}$.\label{Theorem.DownwardLowenheimSkolemIfWellPresented}

\Proof. The idea is to effectively verify the Tarski-Vaught criterion  on the shortest initial elementary cut of the well presented  enumeration of
$\mathcal A$. So, suppose that $\<a_\alpha\st\alpha<\omega_1>$ is the well  presented enumeration of the underlying set of $\mathcal A$. By classical
methods, there is a closed unbounded set of countable initial segments  of this enumeration that form elementary substructures of $\mathcal A$.  Let
$\beta$ be least such that $B=\set{a_\alpha\st\alpha<\beta}$ forms an  elementary substructure ${\mathcal B}\elesub{\mathcal A}$. Thus, $\beta$ is
least such that the set $\set{a_\alpha\st\alpha<\beta}$ satisfies the  Tarski-Vaught criterion in $\mathcal A$. We will argue that $B$ is infinite
time decidable as a set. Given any $a_\xi$, we can generate the sequence  $\<a_\alpha\st\alpha<\xi>$ and for each $\xi'\leq\xi$ we can check whether
$\set{a_\alpha\st\alpha<\xi'}$ satisfies the Tarski-Vaught criterion in  $\mathcal A$. To check this, we use the writable real coding the language to
generate a list of all formulas $\varphi$ in the language. For every  such formula $\varphi$ and every finite sequence
$a_{\alpha_0},\ldots,a_{\alpha_n}$ with each $\alpha_i<\xi'$, we use the  decidability of $\mathcal A$ to inquire whether $\exists
x\,\varphi(x,a_{\alpha_0},\ldots,a_{\alpha_n})$ is true in $\mathcal  A$. If so, then we check that there is some $\alpha<\xi'$ with
$\varphi(a_\alpha,a_{\alpha_0},\ldots,a_{\alpha_n})$ true in $\mathcal  A$. These checks will all be satisfied if and only if
$\set{a_\alpha\st\alpha<\xi'}$ satisfies the Tarski-Vaught criterion.  Consequently, if such a $\xi'$ exists with $\xi'\leq\xi$, then by the
minimality of $\beta$, it must be that $\beta\leq\xi'$, and so $a_\xi$ is not in  $B$. If no such $\xi'$ exists up to $\xi$, then $\xi<\beta$ and so
$a_\xi\in B$. Therefore, as a set, $B$ is decidable. The corresponding model $\mathcal B$ is therefore a decidable model, and a countable elementary
substructure of $\mathcal A$, as desired.\QED

Finally, we have a strong violation to the infinite time computable  downward Lowenheim-Skolem Theorem, when it comes to computable models.  For
infinite time Turing machines, a {\df computation snapshot} is a real coding the  complete description of a machine configuration, namely, the
program  that the machine is running, the head position, the state and the contents of  the cells.

\Theorem. There is an infinite time computable structure with underlying set $\R$ having no infinite time computable proper elementary
substructure.\label{Theorem.ComputableModelOnRwithNoElementarySubstructures}

\Proof. Define the relation $U_p(x,y)$ if $y$ codes the computation  sequence of program $p$ on input $x$ showing it to have been accepted.  That is,
$y$ codes a well-ordered sequence of computation snapshots  $\<y_\alpha\st\alpha\leq\beta>$, such that (i) the first snapshot $y_0$  is the starting
configuration of the computation of program $p$ on input $x$; (ii)  successor snapshots $y_{\alpha+1}$ are updated correctly from the prior  snapshot
$y_\alpha$ and the operation of $p$; (iii) limit snapshots $y_\xi$ correctly show the head on the left-most cell in the {\it limit} state,  with the
tape updated correctly from the prior tape values in  $\<y_\alpha\st\alpha<\xi>$; and lastly, (iv) the final snapshot  $y_\beta$ shows that the
computation halted and accepted the input. This is a computable  property of $\<p,x,y>$, since one can computably verify that $y$ codes  such a well
ordered sequence of snapshots by counting through the underlying order  of $y$ and systematically checking each of the requirements. So the structure
${\mathcal R}=\<\R,U_p>_{p\in\N}$ is a computable structure. (One could  reduce this to a finite language with a trinary predicate $U(p,x,y)$, by
regarding programs as reals and ensuring that the programs are  necessarily in any elementary substructure.)

Suppose that there is a computable proper elementary substructure ${\mathcal A}\elesub{\mathcal R}$. Let $p_0$ be a program deciding the  underlying
set $A$ of $\mathcal A$. Since every real $a\in A$ is accepted by $p_0$, there will be a real $y$ in $\R$ coding the computation sequence and
witnessing $U_{p_0}(a,y)$. Thus, ${\mathcal A}\satisfies\forall a\,\exists  y\,U_{p_0}(a,y)$. By elementarity ${\mathcal A}\elesub{\mathcal R}$, we
conclude that $\mathcal R$ also satisfies this assertion. So every real is accepted  by $p_0$. Thus, $A=\R$ and the substructure is not a proper
substructure after all.\QED

Since this model is only infinite time computable and not infinite time  decidable (the halting problem $0^\jump$ is expressible in the  $\Sigma_1$
diagram), the following question remains open:

\Question. Is there an infinite time decidable model with underlying  set $\R$ having no proper infinite time computable elementary  substructure?

Such a model would be a very strong counterexample to the infinite time  computable downward Lowenheim-Skolem Theorem.

\Section Computable quotient presentations \label{Section.ComputableQuotientPresentations}

We have defined that a structure $\mathcal A$ has an infinite time {\df computable presentation} if it is isomorphic to an infinite time  computable
structure, a structure $\<A,\ldots>$ whose underlying set $A\of\R$ is a decidable set of reals and whose functions and relations are uniformly
computable from their G\"odel codes and their respective input. We define now that a structure $\mathcal A$ has an infinite time computable {\df
quotient presentation} if there is an infinite time computable structure $\mathcal B=\<B,\ldots>$ and an infinite time computable equivalence
relation $\equiv$ on $B$ which is a congruence with respect to the functions and  relations of $\mathcal B$, such that $\mathcal A$ is isomorphic to
the  quotient structure $\mathcal B/{\equiv}$. For example, every computable  structure has a computable quotient presentation, using the equivalence
relation of identity, but there are other more interesting examples. The difference  has to do with the two possibilities in first order logic of
treating  $=$ as a logical symbol, insisting that it be interpreted as identity in a  model, or treating it axiomatically, so that it can be
interpreted  merely as an equivalence relation. The natural question here, of course, is whether  the two notions coincide.

\Question. Does every structure with an infinite time computable quotient presentation have an infinite time computable
presentation?\label{Question.QuotientPresentation}

This is certainly true in the context of finite time computability,  because one can build a computable presentation by using the least  element of
each equivalence class. More generally, for the same reason, it is true in  the infinite time context for structures having a quotient presentation
whose underlying set is contained in the natural numbers. Specifically, if  ${\mathcal A}=\<A,\ldots,\equiv>$ is computable and $A\of\N$, where
$\equiv$ is a congruence on $\mathcal A$, then ${\mathcal A}/{\equiv}$ has a  computable presentation. This is because the function $s$, mapping
every $n\in A$ to the least element $s(n)$ in the equivalence class of $n$, is computable. To  compute $s(n)$, one may simply try out all the smaller
values in turn to discover the least representative. It follows that the set $B=\set{s(n)\st n\in A}$ is a computable choice set for the collection
of equivalence classes. For any relation symbol $R$ in the language of $\mathcal A$,  we may now naturally define $R^{\mathcal B}(\vec n)\iff
R^{\mathcal  A}(\vec n)$; and for any function symbol $f$ we define $f^{\mathcal B}(\vec n)=s(f^{\mathcal A}(\vec n))$. These are clearly computable
functions  and relations, and since $\equiv$ is a congruence, it follows that ${\mathcal A}/{\equiv}$  is isomorphic to $\mathcal B$, as desired.
This argument shows more  generally that if a structure has a computable quotient presentation $\<A,\ldots,\equiv>$, and there is a computable
function $s$ mapping  every element to a representative for its equivalence class, then the quotient structure  ${\mathcal A}/{\equiv}$ has a
computable presentation. (Note: the range  of such a computable choice function will be decidable, because it is precisely  the collection of $x$ in
the original structure for which $s(x)=x$.) Such a function $s$ is like a computable choice function on the equivalence  classes.

In the general infinite time context, of course, one does not expect  necessarily to be able to effectively compute representatives from each
equivalence class. In fact, we will show that the answer to Question  \ref{Question.QuotientPresentation} is independent of \ZFC. In order to
illustrate the ideas, let us begin with the simple example of the uncountable well  order $\<\omega_1,<>$.

\Theorem. \ \begin{enumerate}
 \item The uncountable well-ordered structure $\<\omega_1,<>$ has an  infinite time computable quotient presentation.
 \item It is relatively consistent with \ZFC\ that $\<\omega_1,<>$ has  no infinite time computable presentation.
\end{enumerate}\label{Theorem.ConsistentOmega1HasNoPresentation}

\Proof. For the first claim, observe that the structure $\<\WO,<,{\equiv}>$ is an infinite time computable quotient presentation of $\<\omega_1,<>$.
For any $x\in\WO$, the equivalence class $[x]_\equiv$ is exactly the  set of reals coding the same ordinal as $x$, and so $\<\WO,<>/{\equiv}$  is
isomorphic to $\<\omega_1,<>$, as desired.

For the second claim, observe that by forcing, one may easily collapse  $\omega_1^L$ and add sufficient Cohen generic reals, so that in the  forcing
extension $V[G]$ we have that $\omega_1^L$ is countable and the \CH\  fails. By Lemma \ref{Lemma.NoDecidableSetsSizeOmega1}, therefore, the  model
$V[G]$ has no computable structures of size $\omega_1$. In particular, in  $V[G]$ the structure $\<\omega_1,<>$ has no computable presentation, as
desired.\QED

Thus, it is consistent that the answer to Question \ref{Question.QuotientPresentation} is negative.  We turn now to the  possibility of a positive
answer. Let us begin with a positive answer for the specific structure  $\<\omega_1,<>$.

\Theorem. If\/ $\omega_1=\omega_1^L$ (a consequence of $V=L$), then the  structure $\<\omega_1,<>$ has an infinite time computable
presentation.\label{Theorem.Omega1LHasComputablePresentation}

\Proof. We already observed after Theorem \ref{Theorem.Lomega1HasComputablePresentation} that $\<\omega_1^L,<>$  has a computable presentation using
the $L$-codes for ordinals.\QED

It seems likely that one doesn't really need the failure of \CH\ in the  proof of Theorem \ref{Theorem.ConsistentOmega1HasNoPresentation}, and we
suspect that the particular structure $\<\omega_1,<>$ has a computable presentation if and only if $\omega_1=\omega_1^L$. That is, we suspect  that
the converse of Theorem \ref{Theorem.Omega1LHasComputablePresentation} also  holds.

\Corollary. The question of whether the structure $\<\omega_1,<>$ has  an infinite time computable presentation is independent of \ZFC.

\Proof. On the one hand, by Theorem \ref{Theorem.ConsistentOmega1HasNoPresentation} it is relatively consistent that $\<\omega_1,<>$ has no
computable presentation. On the other hand, if $V=L$ or merely  $\omega_1^L=\omega_1$, then $\<\omega_1,<>$ has a computable  presentation.\QED

Rather than studying just one structure, however, let us now turn to  the possibility of a full positive solution to Question
\ref{Question.QuotientPresentation}. Under $V=L$, one has a full affirmative answer.

\Theorem. If $V=L$, then every structure with an infinite time computable quotient presentation has an infinite time computable
presentation.\label{Theorem.V=LimpliesYes}

\Proof. Assume $V=L$ and suppose that ${\mathcal A}=\<A,\ldots,\equiv>$  is a computable structure, where $\equiv$ is a congruence with respect  to
the rest of the structure. We would like to show that ${\mathcal  A}/{\equiv}$ has a computable presentation. Our argument will be guided  by the
idea of building a computable presentation of ${\mathcal A}/{\equiv}$ by  selecting the $L$-least representatives of each equivalence class. We  will
not, however, be able to do exactly this, because we may not be able to  recognize that a given real is the $L$-least representative of its
equivalence class. Instead, we will attach an escort $y$ to every such $L$-least  representative $x$ of an equivalence class $[x]$, where $y$ codes
an  ordinal sufficiently large to allow us computably to verify that $x$ is the  $L$-least representative of its equivalence class. We will then
build  the computable presentation out of these escorted pairs $\<x,y>$.

First, for simplicity, consider the case that $\mathcal A$ is a relational structure. Let $B$ be the set of pairs $\<x,y>$ such that  $y$ is an
$L$-code for the least ordinal $\alpha$ such that $x$ is an element of $L_\alpha$ and  $L_\alpha$ satisfies that $x$ is in $A$, that ``$\omega_1$
exists'' and that $x$ is the $L$-least real that is equivalent to $x$. The assertions about membership in $A$ or equivalence can be expressed in
$L_\alpha$ using the programs that compute these relations. Note that because $L_\alpha\satisfies``\omega_1$ exists,'' all the computations for reals
in $L_\alpha$ either halt or repeat  before $\alpha$, and so $L_\alpha$ has access to the full, correct computations for the reals in $L_\alpha$.

We claim that $B$ is decidable. First, the set of $L$-codes is decidable. Next, given that $y$ is the $L$-code of an ordinal $\alpha$,  we can by
Theorem \ref{Theorem.CodingLalpha} compute a code for the whole  structure $L_\alpha$, and so questions of satisfaction in this  structure will be
decidable. Next, we can check that $x$ is an element of $L_\alpha$, and  that $L_\alpha$ satisfies all those other properties, as desired. Checking
that $\alpha$ is least with those properties amounts to checking that $L_\alpha$ thinks there is no $\beta$ having an $L$-code that works.

Next, observe that if $\<x,y>\in B$, then $x$ really is the $L$-least  representative of $[x]$ in $\mathcal A$. The reason is that if $z\equiv  x$
and $z$ precedes $x$ in the $L$ order, then $z$ would be in $L_\alpha$  also, where $y$ codes $\alpha$, and so $L_\alpha$ would know that $z$
precedes $x$. And it is correct about whether $z\equiv x$, since it has the  computation checking this. The point is that $L_\alpha$ can see $x$ and
all its $L$ predecessors, and it knows whether they are equivalent or not. So $L_\alpha$ will be correct about whether $x$ is the L-least
representative of $[x]$.

Finally, we put a structure on $B$ as follows. For a relation symbol  $R$, let $R^{\mathcal B}(\<x_0,y_0>,\ldots,\<x_n,y_n>)$ hold if and  only if
$R^{\mathcal A}(x_0,\ldots,x_n)$, which is computable. For each $a\in  A$, there is an $L$-least representative $x$ in $[a]$, and a least  ordinal
$\alpha$ large enough so that $x$ is in $L_\alpha$ and $L_\alpha$  satisfies all those tests. If $y$ is the $L$-code of $\alpha$, then  $\<x,y>$ will
be in $B$. By mapping $[a]$ to $\<x,y>$, it is clear that ${\mathcal  A}/{\equiv}$ is isomorphic to $\mathcal B$, providing a computable
presentation.

When the language has function symbols, we define  $f^{\mathcal B}(\<x_0,y_0>,\ldots,\<x_n,y_n>)=\<x,y>$, where $x$ is the $L$-least member of
$f^{\mathcal A}(x_0,\ldots,x_n)$ and $y$ is the $L$-code for which $\<x,y>\in B$. The point now is that since $f^{\mathcal A}(x_0,\ldots,x_n)$ is the
result of a computation in $L_\alpha$, where  $\alpha$ is the largest of the ordinals arising from $y_0,\ldots,y_n$,  with the structure $L_\alpha$,
we will be able to find the $L$-least member $x$  of the corresponding equivalence class and the $L$-code $y$ putting $\<x,y>$ into $B$. Thus, we
will be able to compute this information from $\<x_0,y_0>,\ldots,\<x_n,y_n>$, and so $f^{\mathcal B}$ is a computable function. Once again ${\mathcal
A}/{\equiv}$ is isomorphic to $\mathcal B$, as desired.\QED

The argument does not fully use the hypothesis that  $V=L$, but rather only that $A\of L$, since in this case we might as  well live inside $L$. In
particular, any structure that has a computable quotient  presentation using only writable reals or even accidentally writable reals, has a
computable presentation.

\Corollary. The answer to Question \ref{Question.QuotientPresentation}  is independent of \ZFC.

\Proof. By Theorem \ref{Theorem.ConsistentOmega1HasNoPresentation}, it  is relatively consistent that there is a structure with a computable quotient
presentation, but no computable presentation. On the other hand, by  Theorem \ref{Theorem.V=LimpliesYes}, it is also relatively consistent that every
structure with a computable quotient presentation has a computable  presentation.\QED

Another way to express what the argument shows is the following. Let us  say that a function $f\from\R\to\R$ is {\df semi-computable} if its graph is
semi-decidable.

\Theorem. If $V=L$ and $\equiv$ is an infinite time computable equivalence relation on a decidable set, then there is a semi-computable function $f$
such that $x\equiv y$ if and only if $f(x)=f(y)$. Succinctly, every  computable equivalence relation on a decidable set reduces to equality  via a
semi-computable function.\label{Theorem.ReducingComputableEquivTo=}

\Proof. Suppose $\equiv$ is an infinite time decidable equivalence relation on $\R^L$. Let $f(u)=\<x,y>$ where $x$ is the $L$-least member of the
equivalence class $[u]_\equiv$ and $y$ is the $L$-code of the least  $\alpha$ such that $x\in L_\alpha\satisfies``\omega_1$ exists.'' The  relation
$f(u)=\<x,y>$ is decidable, since given $u$ and $\<x,y>$, we can computably verify that $u\equiv x$ and that $y$ is the $L$-code of an  ordinal
$\alpha$; if so, we can compute a code for $L_\alpha$, and from this  code we can check whether $\alpha$ is least such that $x\in
L_\alpha\satisfies``\omega_1$ exists''  and $x$ is the $L$-least member of its equivalence class. The structure $L_\alpha$  is correct about this
because it has all the earlier reals in the $L$-order and it has the full computations determining whether they are equivalent to $x$. So $f$ is
semi-computable. Finally, notice that $u\equiv v$ if and  only if $f(u)=f(v)$, since the value of $f$ depended only on the equivalence  classes
$[u]=[v]$.\QED

This observation opens up a number of natural questions for further  analysis. One naturally wants to consider computable reductions, for  example,
rather than semi-computable reductions. What is the nature of the  resulting reducibility hierarchy? To what extent does it share the  features of
the hierarchy of Borel equivalence relations under Borel reducibility? For  starters, can one show that there is no computable reduction of the
relation $E_0$ (eventual equality of two binary strings) to equality?

On a different topic, Theorem \ref{Theorem.V=LimpliesYes} will allow us to show that a positive answer to the following question is consistent with
\ZFC.

\Question. Does every infinite time decidable structure have an infinite time computable presentation? \label{Question.Presentable}

While this question remains open, we offer two partial solutions. First, we show in Theorem \ref{Theorem.ComputablePresentationWithOneFunction} that
when the language is particularly simple, the answer is affirmative. Second, we show in Theorem \ref{Theorem.DecidableImpliesComputablyPresentable}
that a fully general affirmative answer, for all languages, is consistent with \ZFC. We don't know whether a negative answer is consistent with \ZFC.

\Theorem. In a purely relational language, or in a language with only relation symbols plus one unary function symbol, every infinite time decidable
model has an infinite time computable presentation.\label{Theorem.ComputablePresentationWithOneFunction}

\Proof. In a purely relational language, every decidable structure is already computable. So let us suppose that $\cal A$ is an infinite time
decidable structure in a language with relation symbols plus one unary function symbol $f$. We assume that the language is computably presented, so
that $\singleton{\gcode{f\,}}$ is decidable. For each $a\in{\cal A}$, let $a^*$ be the real coding the list
$\<a,\gcode{f\,},f(a),f^2(a),f^3(a),\ldots>$. Let $\cal A^*$ be the set of all such $a^*$. This is an infinite time decidable set, because if we are
given a real $x$ coding $\<x_0,x_1,x_2,\ldots>$, we can check whether $x_0\in\cal A$ using the fact that the underlying set of $\cal A$ is decidable;
we can check whether $x_1=\gcode{f}$ using the decision algorithm for the language, and after this, we can check whether $x_2=f(x_0)$, $x_3=f(x_2)$
and so on, using the decidability of $\cal A$. So we can check whether $x=a^*$ for some $a$. Next, we put a structure on $\cal A^*$. For each
relation symbol $U$ of $\cal A$, define $U$ on $\cal A^*$ by $U(a_1^*,\ldots,a_n^*)$ if and only if $U(a_1,\ldots,a_n)$. This is computable because
$a$ is computable from $a^*$. Next, define $f^{\cal A^*}(a^*)=(f(a))^*=\<f(a),\gcode{f},f^2(a),f^3(a),f^4(a),\ldots>$. The point is that this is
computable from $a^*$, since $a^*$ lists all this information directly. So the structure $\cal A^*$ is computable (and decidable). Since $a\mapsto
a^*$ is clearly an isomorphism, this proves the theorem.\QED

If the language involves countably many unary function symbols and there is a writable real listing the G\"odel codes of these function symbols, then
a similar construction, using $a^*=\oplus\set{\tau(a)\st \tau\hbox{ is a term}}$, would provide a computable presentation. This idea, however, does
not seem to work with binary function symbols.

\Theorem. It is relatively consistent with \ZFC\ that all infinite time decidable  structures are infinite time computably presentable. Thus, it is
consistent with \ZFC\ that the answer to Question \ref{Question.Presentable} is yes.\label{Theorem.DecidableImpliesComputablyPresentable}

\Proof. Suppose that $\cal A$ is an infinite time decidable structure. Augment the language by adding a constant symbol for every element of $\cal
A$, and let $\cal A^*$ be the set of all terms in this expanded language. The function symbols have their obvious interpretations and are computable;
the relations have their natural interpretations and are decidable (since $\cal A$ is decidable). Define $t_1\equiv t_2$ if $\cal A\satisfies
t_1=t_2$. This is a computable equivalence relation, because $\cal A$ is decidable. Since ${\cal A}^*/{\equiv}$ is isomorphic to ${\cal A}$, we have
provided an infinite time computable quotient presentation for $\cal A$. By Theorem \ref{Theorem.V=LimpliesYes}, it is relatively consistent with
\ZFC\ that all such structures have a computable presentation.\QED

We note that in Theorem \ref{Theorem.DecidableImpliesComputablyPresentable}, the computable presentation may involve a computable translation of the
language.

\Section The infinite time analogue of Schr\"{o}der-Cantor-Bernstein-Myhill

In this section, we prove the infinite time computable analogues of the  Schr\"{o}der-Cantor-Bernstein theorem and Myhill's theorem. With the
appropriate hypotheses, as in Theorems \ref{Theorem.CantorSchroederBernstein} and \ref{Theorem.MyhillsTheorem}, the proofs go through with a
classical argument. But let us first discuss the need for careful hypotheses. The  usual proofs of Myhill's theorem and the
Cantor-Schr\"oder-Bernstein theorem involve iteratively applying the functions in a zig-zag pattern between  the two sets. And one of the useful
properties of computable injective functions in the classical finite time context is that their  corresponding inverse functions are also
automatically computable: to compute $f^\inverse(b)$, one simply searches the domain for an $a$ such that  $f(a)=b$. Unfortunately, this method does
not work in the infinite time context, where we generally have no ability to effectively enumerate the domain,  and indeed, there are infinite time
{\df one-way} computable functions  $f$, meaning that $f$ is computable but $f^\inverse$ is not computable. An  easy example of such a function is
provided by the Lost Melody Theorem \ref{Theorem.LostMelodyTheorem}, where we have a real $c$ such that  $\{c\}$ is decidable, but $c$ is not
writable. It follows that the function $c\mapsto 1$ on the singleton domain $\{c\}$ is computable, but its  inverse is not. Building on this, we can
provide a decidable counterexample to a direct infinitary computable analogue of Myhill's theorem.

\Theorem. In the infinite time context, there are decidable sets $A$  and $B$ with computable total injections $f:\R\to\R$ and $g:\R\to\R$  such that
$x\in B\iff f(x)\in A$ and $x\in A\iff g(x)\in B$, but there is no computable bijection  $h:A\to B$.

\Proof. Let $A=\N$ and $B=\N\union\singleton{c}$, where $c$ is the real  of the Lost Melody Theorem. Define $f(c)=0$, $f(n)=n+1$ for $n\in\N$ and
otherwise $f(x)=x$. Clearly, $f$ is a computable total injection and  $x\in B\iff f(x)\in A$. To help define $g$, for any real $x$ (infinite  binary
sequence), let $x^*$ be the real obtained by omitting the first digit,  and let $x^{*(n)}$ be the real obtained by omitting the first $n$  digits.
Now let $g(c)=c^*$ and more generally $g(c^{*(n)})=c^{*(n+1)}$, and otherwise $g(x)=x$. This function $g$ is clearly total and injective,  and it is
computable because given any $x$, we can by adding various finite  binary strings to the front of $x$ determine whether $x=c^{*(n)}$ for  some $n$
and thereby compute $g(x)$. Since $c$ is not periodic, we have  $c\notin\ran(g)$ and $x\in A\iff g(x)\in B$. Finally, there can be no  computable
onto map from $A$ to $B$, since $c$ is not the output of any computable function  with natural number input.\QED

In this example, the function $f\restrict B$ is actually a computable  bijection in the converse direction, from $B$ to $A$, but this doesn't
contradict the theorem because $f^\inverse$ is not computable from $A$ to $B$,  since it maps $0$ to $c$. What we really want in the infinitary
context  is not merely a computable bijection from $A$ to $B$, but rather a computable  bijection whose inverse is also computable, so that the
relation is symmetric. The next example shows that we cannot achieve this even when we have  computable bijections in both directions.

\Theorem. \label{Theorem.bijection} In the infinite time context, there  are decidable sets $A$ and $B$ with computable bijections
$f:A\longrightarrow B$ and $g:B\longrightarrow A$, for which there is no computable bijection  $h:A\longrightarrow B$ whose inverse $h^{-1}$ is also
computable.

To construct $A$ and $B$, we first generalize the Lost Melody Theorem  by recursively building a sequence of reals which can each be  recognized, but
not written, by an infinite time Turing machine using the preceding reals  of the sequence.

\SubLemma. There exists a sequence $\langle d_{k}\mid k\in\omega\rangle$ of reals such that
\begin{enumerate}
 \item  for each $k$, the real $d_k$ is not writable from $\langle  d_i\mid i<k\rangle$ and
 \item  there is an infinite time program which, for any $z$ and any  $k$, can decide
   on input $\langle d_0, d_1, \ldots, d_{k-1}, z\rangle$ whether  $z=d_k$.
\end{enumerate}

\Proof. The repeat-point of a computation is the least ordinal stage by  which the computation either halts or enters a repeating loop from  which it
never emerges.  For each $k\geq 0$, let $\delta_k$ be the supremum of  the repeat-points of all computations of the form  $\varphi_p(\<d_{i}\mid
i<k>)$. Note that $\delta_{k}$ is countable in $L$. Let $\beta_{k}$ be the  smallest ordinal greater than $\delta_k$ such that $L_{\beta_{k}
+1}\models `\beta_{k}$ is countable'. Finally, let $d_k$ be the $L$-least real  coding $\beta_k$. The real $d_k$ is not writable on input $\langle
d_{i}\mid i<k\rangle$, for if it were, then we could solve the halting problem  relative to $\<d_{i}\mid i<k>$ by writing $d_k$ and using it to check
whether any given program halts within $\beta_k$ steps on input $\<d_{i}\mid i<k>$.  Next, on input $\<d_0, d_1, \ldots, d_{k-1}, z>$, let us explain
how to determine whether $z=d_k$. We first check whether $z$ codes an ordinal  $\alpha$, and if so, we simulate every computation
$\varphi_{p}(\<d_{i}\mid i<k>)$ for $\alpha$ many steps. By inspecting these computations, we  can verify that they all halt or repeat by stage
$\alpha$, and thereby  verify that $\alpha\geq \delta_{k}$. By Theorem \ref{Theorem.CodingLalpha}, we can  now write down a real coding
$L_{\alpha+1}$ and verify that $z$ is the $L$-least code for $\alpha$ in $L_{\alpha+1}$. If all these tests are passed,  then $z=d_k$.\QED

\Proof Proof of Theorem \ref{Theorem.bijection}. We use the sequence  $\langle d_k\mid k\in\omega\rangle$ to construct a bi-infinite sequence
$\langle c_k\mid k\in\Z\rangle$ as follows: for $k>0$ let $c_{k}$ be a real  coding $\langle d_{i}\mid i<k\rangle$ in the usual manner, and for
$k\leq 0$, let $c_k=k$. Let $A=\{c_{2k}\mid k\in\Z\}$ and $B=\{c_{2k+1}\mid k\in\Z\}$, and define bijections $f:A\to B$ by $f:c_{2k}\mapsto c_{2k-1}$
and $g:B\to A$ by $f:c_{2k+1}\mapsto c_{2k}$. It follows immediately from the  definition of $c_k$ that $f$ and $g$ are computable.

We next show that $A$ is decidable.  Given a real $z$, we first verify  that either $z$ is an even integer less than or equal to zero, in which  case
we accept it immediately, or else it codes a sequence $\langle z_0, \ldots  z_{n-1}\rangle$ of even length, in which case we use the lemma
iteratively to verify that $z_i=d_i$ for each $i<n$.  Since the real $z$ is an element  of $A$ if and only if it passes this test, $A$ is decidable.
Similarly,  $B$ is decidable.

We conclude by showing that if $h:A\longrightarrow B$ is a bijection,  then $h$ and $h^{-1}$ cannot both be computable. From clause (1) of the  lemma
and the definition of $c_n$, it follows that for positive $n$, $c_n$ cannot  be written by any machine on input $c_k$ if $k<n$.  Thus, if $h$ is
computable, then $h(c_2)$ must equal $c_k$ for some $k<2$.  But then  $h^{-1}(c_k)=c_2$ so $h^{-1}$ is not computable.\QED

\Corollary. In the infinite time context, there are decidable sets $A$  and $B$ and a computable permutation $\pi:\R\to\R$ such that $\pi\image  A=B$
and $\pi\image B=A$, but there is no computable bijection $h:A\to B$ for which  $h^\inverse$ is also computable.

\Proof. Let $A$ and $B$ be as in the proof of Theorem \ref{Theorem.bijection}. Since  $A$ and $B$ are disjoint, the function $\pi=f\union g\union
\id$, where  we use the identity function outside $A\union B$, is a permutation of $\R$.  Since $f$ and $g$ are computable and $A$ and $B$ are
decidable, it follows that $\pi$ is computable. Since $\pi\image A=f\image A=B$ and $\pi\image B=g\image B=A$, the proof is completed by mentioning
that Theorem \ref{Theorem.bijection} shows that there is no computable bijection from $A$ to $B$ whose inverse is also computable.\QED

If one assumes merely that the inverses of the injections are computable, then this is insufficient to get a computable bijection:

\Theorem. In the infinite time context, there are semi-decidable sets  $A$ and $B$ with computable injections $f:A\to B$ and $g:B\to A$ whose
inverses are also computable, such that there is no computable bijection $h:A\to B$.

\Proof. In fact, there will be no computable surjection from $A$ to  $B$. Let $A=\N$ be the set of all natural numbers, and let $B=0^\jump=\set{p\st
\varphi_p(0)\converges}$ be the infinite time halting problem. Define  an injective function $f:A\to B$ by setting $f(n)$ to be the $n^\th$  program
on a decidable list of obviously halting programs (such as the program with  $n$ states and all transitions leading immediately to the {\it halt}
state). The function $f$ is clearly computable, and by design its inverse is also  computable and $\ran(f)$ is decidable. Conversely, construing
programs  as natural numbers, the inclusion map $g:B\to A$ is a computable injection whose  inverse is also computable, since
$\dom(g^\inverse)=0^\jump$ is  semi-decidable. So we have defined the required computable injections. Suppose now that $h:A\to B$ is a computable
surjection of $\N$ to $0^\jump$. In this  case, an infinite time computable function could systematically compute all the  values $h(0)$, $h(1)$, and
so on, and thereby write $0^\jump$ on the  tape. This contradicts the fact that $0^\jump$ is not a writable real. So there  can be no such computable
bijection from $A$ to $B$.\QED

In the classical finite time context, of course, there is a computable  bijection between $\N$ and $0'$ (or any infinite c.e.~set), mapping  each $n$
to the $n^\th$ element appearing in the canonical enumeration of it. This  idea does not work in the infinitary context, however, because the
infinitary halting problem $0^\jump$ is not computably enumerated in order type $\omega$, but rather in the order type $\lambda$ of the clockable
ordinals. And $\lambda$ is not a writable ordinal, so there is no way to effectively produce a  real coding it.

Finally, with the right hypotheses, we prove the positive results, starting with the effective content of the Cantor-Schr\"oder-Bernstein  theorem.

\Theorem. In the infinite time context, suppose that $A$ and $B$ are  semi-decidable sets, with computable injections $f:A\rightarrow B$ and
$g:B\rightarrow A$, whose inverses are computable and whose ranges are decidable. Then there is a computable bijection $h:A\to B$ whose inverse is
computable.\label{Theorem.CantorSchroederBernstein}

\Proof. Let $A_0$ be the set of $a$ such that there is some finite zig-zag pre-image $(g^\inverse f^\inverse)^kg^\inverse(a)\notin\ran(f)$ for
$k\in\N$. Our hypotheses ensure that this set is infinite time decidable, since we can  systematically check all the corresponding pre-images to see
that when and if they stop it was because they landed outside $\ran(f)$ in $B$. The  usual proof of the Cantor-Schr\"oder-Bernstein theorem now shows
that  the function $h=(g^\inverse\restrict A_0)\union(f\restrict A\minus A_0)$ is a bijection between $A$ and $B$. Note that $h$ is computable
because $g^\inverse$ and $f$ are each computable, $A_0$ is decidable, and $A$ is semi-decidable. To see that $h^\inverse$ is computable, let
$B_0=g^\inverse A_0$ and observe that $h^\inverse=(g\restrict B_0)\union(f^\inverse\restrict B\minus B_0)$. Since these components  are each
computable, $h^\inverse$ is computable and the proof is  complete.\QED

We may drop the assumption that $A$ and $B$ are semi-decidable if we  make the move to total functions, as in the classical Myhill's theorem.  Define
that a set of reals $A$ is {\df reducible} to another set $B$ by the  function $f:\R\to\R$ if $x\in A\iff f(x)\in B$.

\Theorem. In the infinite time context, suppose that $A$ and $B$ are  reducible to each other by computable one-to-one total functions $f$  and $g$,
whose inverses are computable and whose ranges are decidable. Then there is a computable permutation $\pi:\R\to\R$ with $\pi^\inverse$ also
computable and $\pi\image A=B$.\label{Theorem.MyhillsTheorem}

\Proof. As in Theorem \ref{Theorem.CantorSchroederBernstein}, let $A_0$  be the set of $a$ such that some finite zig-zag pre-image $(g^\inverse
f^\inverse)^kg^\inverse(a)\notin\ran(f)$ for $k\in\N$, and again this  is decidable. Let $\pi=(g^\inverse\restrict  A_0)\union(f\restrict\R\minus
A_0)$. The usual Cantor-Schr\"oder-Bernstein argument shows that this is a  permutation of $\R$. As above, both $\pi$ and $\pi^\inverse$ are
computable. Finally, we have both $x\in A\iff\pi(x)\in B$ and $x\in B\iff \pi(x)\in  A$, since $\pi(x)$ is either $f(x)$ or $g^\inverse(x)$, both of
which  have the desired properties. It follows that $\pi\image A=B$ and the proof is  complete.\QED

%
%
%
%
%
\Section Some infinite time computable transitive models of set theory

Because the power of the machines are intimately connected with well-orders and countable ordinals, it is not surprising that there are  many
interesting models of a set theoretic nature. We have already seen that  the hereditarily countable sets have an infinite time computable  quotient
presentation $\<\HC,{\in},\equiv>$. In addition, we have provided  infinite time computable presentations of the model $\<L_{\omega_1^L},{\in}>$ and
of $\<L_\alpha,{\in}>$, given a real coding $\alpha$. We will now show,  however, that depending on the set theoretic background, one can  transcend
these, by actually producing infinite time decidable presentations of {\it transitive} models of \ZFC, or even \ZFC\ plus large cardinals.

\Theorem. If there is a transitive model of \ZFC, then the smallest  transitive model of \ZFC\ has an infinite time decidable computable
presentation.\label{Theorem.MinimalZFCModel}

\Proof. If there is a transitive model of \ZFC, then there is one satisfying $V=L$. A L\"owenheim-Skolem argument, followed by the Mostowski
collapse, shows that there must be a countable such model, and any such model  will be $L_\alpha$ for some countable ordinal $\alpha$. By minimizing
$\alpha$, we see that there is a smallest transitive model $L_\alpha\satisfies\ZFC$.  Let $c$ be the $L$-code for this minimal $\alpha$. Note that
$\set{c}$  is decidable, since on input $x$, we can check whether it is an $L$-code  for an ordinal $\xi$ such that $L_\xi\satisfies\ZFC$, and if so,
whether $\xi$ is the smallest such ordinal. If so, it must be that $\xi=\alpha$ and  $x=c$. From $c$, we may compute a relation $E$ on $\N$ such that
$\<L_\alpha,{\in}>\iso\<\N,E>$. Let $M$ be the collection of pairs  $\<c,n>$, where $n\in\N$. This is a decidable set, because $\set{c}$ is
decidable. The idea is that $\<c,n>$ represents the set coded by $n$ with respect  to $E$. Define $\<c,n>\mathrel{\bar E}\<c,m>$ if $n\mathrel{E}m$.
This  is a computable relation, because $E$ is $c$-computable and $\singleton{c}$  is decidable. Clearly, $\<M,\bar E>$ is isomorphic to $\<\N,E>$,
which  is isomorphic to $\<L_\alpha,{\in}>$. So we have a computable presentation  of $\<L_\alpha,{\in}>$.

Let us point out that this presentation is nearly decidable, in that we  can decide ${\cal M}\satisfies\varphi[x_1,\ldots,x_n]$ on input
$\gcode{\varphi},x_1,\ldots,x_n$, provided $n\geq 1$. Specifically,  from $x_1$ we can compute the real $c$, and from $c$ we can effectively
enumerate the whole structure $\<M,E>$. Having done so, we can compute whether  ${\cal M}\satisfies\varphi[x_1,\ldots,x_n]$ according to the Tarskian
definition of truth. This method makes fundamental use of the information $c$ that  is present in any of the parameters, so it does not help us to
decide  whether a given sentence holds in ${\cal M}$, if we are not given such a parameter.

To make the model fully decidable, therefore, we assume $\gcode{\in}=c$. Since $\singleton{c}$ is decidable, this language remains decidable
(although no longer enumerable in any nice sense). The point now is that if we  are given a sentence $\sigma$, and the symbol $\in$ appears in it,
then  we can compute the real $c$ from $\gcode{\in}$ and thereby once again effectively enumerate the whole structure ${\cal M}$, allowing us to
compute whether $\sigma$ holds. If $\in$ does not occur in $\sigma$, then $\sigma$ is  an assertion in the language of equality, which either holds
or fails  in all infinite models, and we can computably determine this in $\omega+1$  many steps.\QED

If one allows $\gcode{\in}=c$, then one can actually take the underlying set of $\cal M$ to be $\N$, since if one has already coded  $c$ into the
language, there is no additional need to code $c$ into the individual  elements of $\cal M$. In this case, one has a decidable presentation of  the
form $\<\N,E>$. We caution in this case that the relation $E$ is not  computable, but only computable relative to $c$. This does not prevent the
model from being a computable model, however, since in order to be a computable  model, the relations need only be computable from their G\"odel
codes.  This may be considered to be a quirk in the definition of computable model, but in  order to allow for uncountable languages, we cannot
insist that the  relations of a computable model are individually computable, but rather only computable from their G\"odel codes.

Similar arguments establish:

\Theorem. If there is a transitive model of \ZFC\ with an inaccessible  cardinal (or a Mahlo cardinal or $\omega^2$ many weakly compact  cardinals,
etc.), then the smallest such model has an infinite time decidable  presentation.

\Proof. If there is a transitive model of \ZFC\ plus any of these large cardinal hypotheses, then there is one satisfying $V=L$. Hence, as argued in
Theorem 48, the theory holds in some countable $L_\alpha$. By using the  $L$-code $c$ of the minimal such model, we can build a decidable
presentation as above.\QED

If one wants to consider set theoretic theories inconsistent with $V=L$, then a bit more care is needed.

\Theorem. If there is a transitive model of \ZFC, then there is transitive a model  of $\ZFC+\neg\CH$ with an infinite time decidable computable
presentation.

\Proof. Let $L_\alpha$ be the minimal transitive model of \ZFC. This is  a countable transitive model, and so there is an $L$-least set $G$ in  $L$
such that $G$ is $L_\alpha$-generic for the forcing  $\Add(\omega,\omega_2)^{L_\alpha}$. Thus, $L_\alpha[G]\satisfies\ZFC+\neg\CH$. The set $G$
appears in some countable $L_\beta$, where $\alpha<\beta<\omega_1$. Let $d$ be the $L$-code of the pair $\<\alpha,\beta>$. Thus, $\set{d}$ is
decidable,  because given any real $z$, we can check if $z$ is an $L$-code for a pair  $\<\alpha',\beta'>$ such that $L_{\alpha'}$ is the smallest
model of  \ZFC\ and $\beta'$ is smallest such that $L_{\beta'}$ has an  $L_{\alpha'}$-generic filter $G$ for $\Add(\omega,\omega_2)^{L_{\alpha'}}$.
Using the real $d$, we can compute a relation $E$ on $\N$ such that $\<L_\alpha[G],{\in}>\iso\<\N,E>$. Let $N$ be the set of pairs $\<d,n>$  where
$n\in\N$, and define $\<d,n>\mathrel{\bar E}\<d,m>$ if $n\mathrel{E}m$. Again, this  structure is computable, and it is isomorphic to
$\<L_\alpha[G],{\in}>$, as desired. By taking $\gcode{\in}=d$, the model is decidable as in Theorem  \ref{Theorem.MinimalZFCModel}.\QED

Clearly this method is very flexible; it provides decidable presentations of transitive models of any theory having a transitive  model in $L$.
Nevertheless, we admit that all these examples are a bit strange,  because of their stealthy manner of coding information into the  individual
elements of the underlying set or into the language. So we close this section by  proving that every decidable presentation of a model of \ZFC\ must
use these stealthy measures.

\Theorem. There is no infinite time computable presentation of a transitive model of \ZFC\ with underlying set $\N$ and G\"odel codes of the language
entirely in $\N$.\label{Theorem.NoNiceComputableZFCModels}

\Proof. The operation of an infinite time Turing machine is absolute to  any transitive model of \ZFC\ containing the input. Thus, all  transitive
models of \ZFC\ agree on the elements of the halting problem $0^\jump$. If  ${\cal M}=\<\N,E>$ is a computable presentation of such a model, then
there is some natural number $k$ representing $0^\jump$ in $\cal M$. Assuming that $\gcode{\in}$ is writable, then we can computably determine for
each  natural number $p$ the element $k_p$ representing it in $\cal M$. In this case,  we could compute $0^\jump=\set{p\st k_p\mathrel{E}k}$,
contradicting  the fact that $0^\jump$ is not computable.\QED

Of course, the argument uses much less than \ZFC. It shows that there  can be no computable presentation, using underlying set $\N$ and  writable
presentation of the language, of a transitive model computing $0^\jump$  correctly. For example, it would be enough if the model satisfied
``$\omega_1$ exists,'' or even less, that every infinite time Turing computation  either halted or reached its repeat point. So if one wants
decidable or even computable transitive models of set theory, one must code information  into the elements of the model or the language.

\Section Future Directions

We close this paper by mentioning a number of topics for future research.

{\it Infinitary languages $L_{\omega_1,\omega}$}. In the context of  infinite time computable model theory, it is very natural to consider infinitary
languages, which are still easily coded into the reals. With any  writable structure or for a structure whose domain we can search, one can still
compute the Tarskian satisfaction relation. What other examples and  phenomenon exist here?

{\it Infinite time computable equivalence relation theory}. The idea is  to investigate the analogue of the theory of Borel equivalence  relations
under Borel reducibility. Here, one wants to consider infinite time  computable reductions. Some of these issues are present already in our  analysis
of the computable quotient presentation problem in Section \ref{Section.ComputableQuotientPresentations} and particularly Theorem
\ref{Theorem.ReducingComputableEquivTo=}. How much of the structure of  Borel equivalence relations translates to the infinite time computable
context?

{\it Infinite time computable cardinalities}. The computable cardinalities are the equivalence classes of the decidable sets by the  computable
equinumerousity relation. What is the structure of the computable  cardinalities?

{\it Infinite time computable Lowenheim-Skolem Theorems}. While Theorem  \ref{Theorem.UpwardLowenheimSkolemIfV=L} shows that the infinite time
computable upward Lowenheim-Skolem theorem holds in $L$, our analysis  leaves open the question of whether it is consistent with \ZFC\ that  there
could be a decidable countable model having no size continuum decidable  elementary extension. If so, the infinite time computable upward
Lowenheim-Skolem theorem will be independent of \ZFC. In addition, our analysis does not fully settle the infinite time computable downward
Lowenheim-Skolem  theorem.

\bibliographystyle{alpha}
\bibliography{MathBiblio,HamkinsBiblio}

\end{document}